\documentclass{amsart}
\newtheorem{theorem}{Theorem}[section]

\newtheorem{lemma}[theorem]{Lemma}
\newtheorem{proposition}[theorem]{Proposition}
\newtheorem{algorithm}[theorem]{Algorithm}

\newtheorem*{claim}{Claim}
\numberwithin{equation}{section}
\newcommand{\N}{\mathbb N}
\newcommand{\Z}{\mathbb Z}
\newcommand{\Q}{\mathbb Q}
\newcommand{\R}{\mathbb R}
\newcommand{\C}{\mathbb C}

\newcommand{\A}{\mathbb A}

\newcommand{\slp} {straight--line program}
\newcommand{\slps} {straight--line programs}
\newcommand{\x}{\mathbf{x}}

\def\ifm#1#2{\relax \ifmmode#1\else#2\fi}

\newcommand{\klk}    {\ifm {,\ldots,} {$,\ldots,$}}

\newcommand{\plp}    {\ifm {+\cdots+} {$+\ldots+$}}

\newcommand{\om}[2]   {{#1}_1 \klk {#1}_{#2}}

\newcommand{\spar} {\vskip 0.25cm}

\newcommand{\xon}    {\ifm {\om X n} {$\om X n$}}

\begin{document}

\title[Deformation techniques for sparse systems]{Deformation techniques for sparse systems}%

\author[G. Jeronimo]{Gabriela Jeronimo${}^{1,2}$}%
\author[G. Matera]{Guillermo Matera${}^{2,3}$}
\author[P. Solern\'o]{Pablo Solern\'o${}^{1,2}$}
\author[A. Waissbein]{Ariel Waissbein${}^{4,5}$}

\address{${}^{1}$Departamento de Matem\'atica, Facultad de
Ciencias Exactas y Na\-turales, Universidad de Buenos Aires,
Ciudad Universitaria, Pabell\'on I (1428) Buenos Aires,
Argentina.} \email{jeronimo@dm.uba.ar, psolerno@dm.uba.ar}

\address{${}^{2}$ National Council of Science and Technology (CONICET),
Ar\-gentina.}

\address{${}^{3}$Instituto de Desarrollo Humano,
Universidad Nacional de Gene\-ral Sarmiento, J.M. Guti\'errez 1150
(1613) Los Polvorines, Buenos Aires, Argentina.}
\email{gmatera@ungs.edu.ar}

\address{${}^{4}$CoreLabs, CORE ST,
Humboldt 1967 (C1414CTU) Ciudad de Buenos Aires, Argentina.}
\address{${}^{5}$Doctorado en Ingenier\'\i a, Instituto Tecnol\'ogico
de Buenos Aires, Av. Eduardo Madero 399 (C1106ACD) Ciudad de
Buenos Aires, Argentina.}\email{Ariel.Waissbein@corest.com}

\thanks{Research was partially supported by the
following grants: UBACyT X112 (2004--2007), UBACyT X847
(2006--2009), PIP CONICET 2461, PIP CONICET 5852/05, UNGS 30/3005
and MTM2004-01167 (2004--2007).}



\date{\today}%
\begin{abstract}
We exhibit a probabilistic symbolic algorithm for solving
zero--dimensional sparse systems. Our algorithm combines a
symbolic homotopy procedure, based on a flat deformation of a
certain morphism of affine varieties, with the polyhedral
deformation of Huber and Sturmfels. The complexity of our
algorithm is quadratic in the size of the combinatorial structure
of the input system. This size is mainly represented by the mixed
volume of Newton polytopes of the input polynomials and an
arithmetic analogue of the mixed volume associated to the
deformations under consideration.
\end{abstract}
\maketitle
%
%
\section{Introduction}
Numeric and symbolic methods for computing all solutions of a
given zero-dimensional polynomial system usually rely on
deformation techniques, based on a perturbation of the original
system and a subsequent (numeric or symbolic) path-following
method (see, e.g., \cite{AlGe90}, \cite{BaMu93},
\cite{BlCuShSm98}, \cite{HeKrPuSaWa00}, \cite{Li97},
\cite{Morgan87}). The complexity of such algorithms is usually
determined by geometric invariants associated to the family of
systems under consideration (see, e.g., \cite{CaEm95},
\cite{HuSt95}, \cite{VeGaCo96}, \cite{Rojas99},
\cite{HeKrPuSaWa00}, \cite{JeKrSaSo04}, \cite{CaGiHeMaPa03},
\cite{Schost03}, \cite{Lecerf03}, \cite{PaSa04}), typically in the
form of a suitable (arithmetic or geometric) B\'ezout number (see
\cite{MoSoWa95}, \cite{HuSt95}, \cite{LiWa96}, \cite{RoWa96},
\cite{HuSt97}, \cite{HeJeSaSaSo02}, \cite{PhSo06}).

Sparse elimination theory is concerned with finding bounds for
such B\'ezout numbers in the case of a sparse polynomial system.
Its origins can be traced back to the results by D.N. Bernstein,
A.G. Kushnirenko and A.G. Khovanski (\cite{Bernstein75},
\cite{Kushnirenko76}, \cite{Khovanski78}) that bound the number of
solutions of a polynomial system in terms of certain combinatorial
invariants. More precisely, the Bernstein--Kushnirenko--Khovanski
(BKK for short) theorem asserts that the number of isolated
solutions in the $n$--dimensional complex torus $(\C^*)^n$ of a
polynomial system of $n$ equations in $n$ unknowns is bounded by
the mixed volume of the family of Newton polytopes of the
corresponding polynomials.

Numeric (homotopy continuation) methods for sparse systems are
typically based on a family of deformations called polyhedral
homotopies (\cite{HuSt95}, \cite{VeVeCo94}, \cite{VeGaCo96}).
Polyhedral homotopies preserve the Newton polytope of the input
polynomials and yield an effective version of the BKK theorem (see
e.g. \cite{HuSt95}, \cite{HuSt97}). More precisely, suppose that
we are given a zero--dimensional $(\Delta_1\klk\Delta_n)$--sparse
system defined by $n$ polynomials $f_1\klk f_n$ in $n$ variables,
where $\Delta_1\klk\Delta_n$ are the supports of $f_1\klk f_n$,
and let $V\subset(\C^*)^n$ be the variety defined by the common
zeros of $f_1\klk f_n$ over $(\C^*)^n$. Then a polyhedral homotopy
consists in an algebraic curve $W\subset(\C^*)^{n+1}$ such that
the   projection $\pi:W\to \C^*$ onto the first coordinate is
dominant with generically finite fibers whose degree is the mixed
volume $MV\big(conv(\Delta_1)\klk conv(\Delta_n)\big)$ of the
convex hulls of $\Delta_1\klk\Delta_n$, the identity
$\pi^{-1}(1)=\{1\}\times V$ holds and the first terms of the
Puiseux expansions of the branches of $W$ lying above 0 can be
easily computed. Numerical continuation methods compute the first
terms of these Puiseux expansions and then follow the branches of
$W$ along the interval $[0,1]$ to obtain approximations to all the
points of the input variety $V$.

{From} the symbolic point of view, a family of homotopy algorithms
is based on a flat deformation of a certain morphism of affine
varieties. This deformation, implicitly considered in the papers
\cite{GiHeMoMoPa98}, \cite{GiHaHeMoMoPa97}, is isolated in
\cite{HeKrPuSaWa00} and refined in \cite{Schost03},
\cite{HeJeSaSaSo02}, \cite{BoMaWaWa04}, \cite{PaSa04}, in order to
solve particular instances of a parametric system with a finite
generically--unramified linear projection of "low" degree. More
precisely, let $V$ be a zero--dimensional variety of $\C^n$,
defined by a ``square'' system $f_1=\dots=f_n=0$, and let be given
an algebraic curve $W\subset\C^{n+1}$ and a (dominant) projection
mapping $\pi:W\to\C$ which represent a deformation of $V$. Then,
from a complete description of a generic fiber of the projection
$\pi:W\to\C$, it is possible to compute a complete description of
the input fiber, say $\pi^{-1}(1)=\{1\}\times V$. The complexity
of this procedure can be roughly estimated by the product of two
geometric invariants: the degree of the morphism $\pi$ and the
degree of the curve $W$. The algorithm is nearly optimal in worst
case \cite{CaGiHeMaPa03}, and has good performance over certain
{well--posed} families of polynomial systems of practical interest
(see \cite{HeKrPuSaWa00}, \cite{Schost03}, \cite{BoMaWaWa04},
\cite{CaMaWa06}).

In this article we combine these symbolic techniques, particularly
in the version of \cite{BoMaWaWa04}, with the polyhedral
deformation \cite{HuSt95}, in order to derive a symbolic
probabilistic algorithm for solving sparse zero--dimensional
polynomial systems with quadratic complexity in the size of the
combinatorial structure of the input system. More precisely,
suppose that we are given polynomials $f_1,\ldots,f_n$ of
$\Q[X_1,\ldots,X_n]$ such that the system $f_1=0,\ldots,f_n=0$
defines a zero--dimensional affine subvariety $V$ of $\C^n$.
Denote by $\Delta_1,\ldots,\Delta_n\subset\Z_{\ge 0}^n$ the
supports of $f_1,\ldots,f_n$, and assume that  $0\in\Delta_i$ for
$1\le i\le n$ and the mixed volume $D$ of the Newton polytopes
$conv(\Delta_1),\ldots,conv(\Delta_n)$ is nonzero. Then, given a
``sufficiently generic'' lifting function $\omega$ for
$(\Delta_1,\ldots,\Delta_n)$ and the corresponding mixed
subdivision, we exhibit an algorithm which computes a complete
description of the solution set $V$ of the input system $f_1=0\klk
f_n=0$.

The polyhedral deformation under consideration requires that the
coefficients of the input polynomials satisfy certain generic
conditions (see Section \ref{subsect: polyhedral deform}). For
this reason, we introduce some auxiliary {\em generic} polynomials
$g_1\klk g_n$ with the same supports $\Delta_1\klk\Delta_n$ and
consider the perturbed polynomial system $h_i:= f_i +g_i$ for
$1\le i \le n$. We first solve the system $h_1=0, \dots, h_n = 0$
and then recover the solutions to the input system $f_1 = 0,
\dots, f_n = 0$ by considering a homotopy of type $f_1 + (1-T)
g_1, \dots, f_n + (1-T) g_n$ (in Section \ref{sect: solution
original system}).

The system $h_1 = 0, \dots, h_n = 0$ is solved by considering the
polyhedral homotopy of \cite{HuSt95}. This homotopy introduces a
new variable $T$ and deforms the polynomial $h_i$ by multiplying
each nonzero monomial of $h_i$ by the power of $T$ determined by
the given lifting function $\omega$. From the genericity of the
coefficients of $h_1\klk h_n$ we conclude that the roots of the
resulting parametric system $\widehat h_1 = 0, \dots, \widehat h_n
= 0$ are algebraic functions of the parameter $T$ whose expansions
as Puiseux series can be obtained by ``lifting'' the solutions to
certain zero--dimensional polynomial systems
$h_{1\!,\gamma}^{(0)}=\dots=h_{n\!,\gamma}^{(0)}=0$ associated to
the lower facets $\widehat{C}_\gamma$ of the lifted polytopes
$conv(\widehat{\Delta}_1)\klk conv(\widehat{\Delta}_n)$ defined by
$\omega$. These polynomial systems
$h_{1\!,\gamma}^{(0)}=\dots=h_{n\!,\gamma}^{(0)}=0$ can be easily
solved due to their specific structure, which enables us to use
their solutions as a starting point for our computations (see
Section \ref{subsect: geo sol V0 gamma} for details).

The complexity of our algorithm is mainly expressed in terms of
two quantities related to the combinatorial structure of the input
system: the mixed volume $D:=M\big(conv(\Delta_1)\klk
conv(\Delta_n)\big)$ and certain (nonarchimedean) heights
$E,E^\prime$ associated to our polyhedral deformations. These
heights, which are an arithmetic analogue of the mixed volume $D$
(see \cite{PhSo03}, \cite{PhSo05}), can be bounded in terms of
certain mixed volumes associated to the polyhedral deformation
under consideration, with equality for a generic choice of the
coefficients of the polynomials $\widehat{h}_{i}$ (see Lemma
\ref{lemma: estimate sparse height} below; compare also with
\cite[Theorem 1.1]{PhSo06}). Therefore, we may paraphrase our
complexity estimate as saying that it is {\em quadratic} in the
combinatorial structure of the input system, with a geometric and
an arithmetic component. More precisely, our algorithm requires
$\mathcal{L}n^{O(1)}D\max\{E,E^\prime\}$ arithmetic operations
over $\Q$ (up to polylogarithmic terms), where $\mathcal{L}$ is
the number of arithmetic operations required to evaluate the
polynomials $\widehat{h}_i$ and $f_i+g_i$, and $E$ and $E^\prime$
denote the height of the varieties defined by $\widehat{h}_1=0\klk
\widehat{h}_n=0$ and $f_1+(1-T)g_1=0\klk f_n+(1-T)g_n=0$
respectively.

This improves and refines the estimate of \cite{BoMaWaWa04} in the
case of a sparse system, which is expressed as a fourth power of
$D$ and the maximum of the degrees of the varieties
$\widehat{h}_1=0\klk \widehat{h}_n=0$ and $f_1+(1-T)g_1=0\klk
f_n+(1-T)g_n=0$. We observe that this maximum is an upper bound
for the heights $E$ and $E^\prime$ respectively. On the other
hand, it also improves \cite{Rojas99}, \cite{Rojas00b}, which
solve a sparse system with a complexity which is roughly quartic
in the size of the combinatorial structure of the input system.
Finally, we provide an explicit estimate of the error probability
of all the steps of our algorithm. This might be seen as a further
contribution to the symbolic stage of the probabilistic
seminumeric method of \cite{HuSt95}, which lacks such analysis of
the error probability.
%
%
\section{Notions and notations}\label{sect: notions and notations}
%
%
\subsection{Sparse Elimination}\label{subsect: sparse elimination}
Here we introduce some notions and notations of  convex geometry
and sparse elimination theory (see e.g. \cite{GeKaZe94},
\cite{HuSt95}) that will be used in the sequel.

Let $X_1,\ldots,X_n$ be indeterminates over $\Q$ and write
$X:=(X_1,\ldots,X_n)$. For $q:=(q_1,\ldots,q_n)\in\Z^n$, we use
the notation $X^q:=X_1^{q_1}\cdots X_n^{q_n}$. Let $f:=\sum_{q}
c_q X^q$ be a Laurent polynomial in $\Q[X,
X^{-1}]:=\Q[X_1,X_1^{-1},\ldots,X_n, X_n^{-1}]$. By the {\sf
support} of $f$ we understand the subset of $\Z^n$ defined by the
elements $q\in\Z^n$ for which $c_q\neq 0$ holds. The {\sf Newton
polytope} of $f$ is the convex hull of the support of $f$ in
$\R^n$.

A {\sf sparse polynomial system} with respect to {\em a priori}
fixed finite subsets $\Delta_1,\ldots,\Delta_n$ of $(\Z_{\ge
0})^n$ is defined by polynomials
$$
f_i(X):=\sum_{q\in\Delta_i}a_{i,q}\, X^q\quad (1\le i\le n),
$$
with $a_{i,q}\in\C\setminus\{0\}$ for each $q\in\Delta_i$ and
$1\le i\le n$.

For a finite subset $\Delta$ of $\Z^n$, we denote by
$Q:=conv(\Delta)$ its convex hull in $\R^n$. The usual Euclidean
volume of a polytope $Q$ in $\R^n$ will be denoted by
$vol_{\R^n}(Q)$.

Let $Q_1,\ldots,Q_n$ be convex polytopes in $\R^n$. For
$\lambda_1,\ldots,\lambda_n\in\R_{\ge 0}$, we use the notation
$\lambda_1Q_1+\cdots+\lambda_nQ_n$ to refer to the Minkowski sum
$\lambda_1Q_1+\cdots+\lambda_nQ_n:=\{x\in\R^n: \exists x_1\in
Q_1,\ldots,\exists x_n\in Q_n \mbox{ such that }
x=\lambda_1x_1+\cdots+\lambda_nx_n\}$. Consider the real--valued
function $(\lambda_1,\ldots,\lambda_n)\mapsto
vol_{\R^n}(\lambda_1Q_1+\cdots+\lambda_nQ_n)$. This is a
homogeneous polynomial function of degree $n$ in the $\lambda_i$
(see e.g. \cite[Chapter 7, Proposition \S 4.4.9]{CoLiOS98}). The
{\sf mixed volume} $MV(Q_1,\ldots,Q_n)$ of $Q_1,\ldots,Q_n$ is
defined as the coefficient of the monomial
$\lambda_1\cdots\lambda_n$ in
$vol_{\R^n}(\lambda_1Q_1+\cdots+\lambda_nQ_n)$.

For $i=1,\dots, n$, let $\Delta_i$ be a finite subset of $\Z_{\ge
0}^n$ and let $Q_i:=conv(\Delta_i)$ denote the corresponding
polytope. Let $f_1,\ldots,f_n$ be a sparse polynomial system with
respect to $\Delta_1,\ldots,\Delta_n$. The BKK theorem
(\cite{Bernstein75}, \cite{Kushnirenko76}, \cite{Khovanski78})
asserts that the system $f_1=0,\ldots,f_n=0$ has at most
$MV(Q_1,\ldots,Q_n)$ isolated common solutions in the
$n$--dimensional torus $(\C^*)^n$, with equality for generic
choices of the coefficients of $f_1,\ldots,f_n$. Furthermore, if
the condition $0\in Q_i$ holds for $1\le i\le n$, then
$MV(Q_1,\ldots,Q_n)$ bounds the number of solutions in the
$n$--dimensional affine complex space $\A^n:=\A^n(\C)$ (see
\cite{LiWa96}).

Assume that the union of the sets $\Delta_1,\ldots,\Delta_n$
affinely generate $\Z^n$, and consider the partition of
$\Delta_1,\ldots,\Delta_n$ defined by the relation
$\Delta_i\sim\Delta_j$ if and only if $\Delta_i=\Delta_j$. Let
$s\in\N$ denote the number of classes in this partition, and let
$\mathcal{A}^{(1)},\ldots,\mathcal{A}^{(s)}\subset\Z^n$ denote a
member in each class. Write
$\mathcal{A}:=(\mathcal{A}^{(1)},\ldots,\mathcal{A}^{(s)})$. For
$\ell=1,\dots, s$, let $k_\ell:=\#\{ i : \Delta_i =
\mathcal{A}^{(\ell)}\}$. Without loss of generality, we will
assume that $\Delta_1 =\dots = \Delta_{k_1} = \mathcal{A}^{(1)}$,
$\Delta_{k_1+1} =\dots = \Delta_{k_1+k_2} = \mathcal{A}^{(2)}$ and
so on.

A {\sf cell} of $\mathcal{A}$ is a tuple
$C=(C^{(1)},\ldots,C^{(s)})$ with $C^{(\ell)} \ne \emptyset$ and
$C^{(\ell)}\subset\mathcal{A}^{(\ell)}$ for $1\le\ell\le s$. We
define
\begin{align*}
type(C):=\,&(\dim(conv(C^{(1)})),\ldots,\dim(conv(C^{(s)}))), \\
conv(C):=\,& conv(C^{(1)}+\cdots+C^{(s)}),\\ \#(C)  :=\,&
\#(C^{(1)})+\cdots+\#(C^{(s)}), \\
vol_{\R^n}(C):=\, & vol_{\R^n}(conv(C)).
\end{align*}
A {\sf face} of a cell $C$ is a cell
$\mathcal{C}=(\mathcal{C}^{(1)},\ldots,\mathcal{C}^{(s)})$ of $C$
with $\mathcal{C}^{(\ell)}\subset C^{(\ell)}$ for $1\le \ell\le s$
such that there exists a linear functional $\gamma:\R^n\to \R$
that takes its minimum over $C^{(\ell)}$ at $\mathcal{C}^{(\ell)}$
for $1\le \ell \le s$. One such functional $\gamma$ is called an
{\sf inner normal} of $C$.

A {\sf mixed subdivision} of $\mathcal{A}$ is a collection of
cells $\mathfrak{C }= \{C_1,\ldots,C_m\}$ of $\mathcal{A}$
satisfying conditions (1)--(4) below:
\begin{enumerate}
\item $\dim(conv(C_j)) = n$ for $1\le j\le m$,
\item the intersection $conv(C_i)\cap conv(C_j)\subset\R^n$ is
either the empty set or a face of both $conv(C_i)$ and $conv(C_j)$
for $1\le i<j\le m$,
\item $\bigcup_{j=1}^m conv(C_j)= conv(\mathcal{A})$,
\item $\sum_{\ell=1}^s\dim(conv(C_j^{(\ell)}))=n$ for $1\le j\le
m$.
\end{enumerate}
If $\mathfrak{C}$ also satisfies the condition
\begin{enumerate}
\setcounter{enumi}{4} \item 
$\#(C_j)=n+s$ for $1\le j\le m$,
\end{enumerate}
we say that $\mathfrak{C}$ is a {\sf fine--mixed subdivision} of
$\mathcal{A}$. Observe that, as a consequence of conditions (4)
and (5), for each cell $C_j = (C_j^{(1)}, \dots, C_j^{(s)})$ in a
fine--mixed subdivision the identity $\dim(conv(C_j^{(\ell)})) =
\# C_j^{(\ell)} - 1$ holds for $1\le \ell \le s$.

We point out that a mixed subdivision $\mathfrak{C}$ of
$\mathcal{A}$ enables us to compute the mixed volume of the family
$Q_1=conv(\Delta_1), \dots, Q_n= conv(\Delta_n)$ by means of the
following identity (see \cite[Theorem 2.4.]{HuSt95}):
\begin{equation}\label{eq: mixed volume and volumes}
MV(Q_1, \dots, Q_n) = \sum_{C_i \in \mathfrak{C} \atop \text{type}(C_i) =
(k_1,\dots, k_s)} k_1 ! \dots k_s! \cdot vol_{\R^n}(C_i).
\end{equation}

A fine--mixed subdivision of $\mathcal{A}$ can be obtained by
means of a lifting process as explained in what follows. For $1\le
\ell \le s$, let $\omega_\ell: \mathcal{A}^{(\ell)}\to\R$ be an
arbitrary function. The tuple $\omega:=(\omega_1,\ldots,\omega_s)$
is called a {\sf lifting function} for $\mathcal{A}$. Once a
lifting function $\omega$ is fixed, the graph of any subset
$C^{(\ell)}$ of $\mathcal{A}^{(\ell)}$ will be denoted by
$\widehat{ C}^{(\ell)}:=\{(q,\omega_\ell(q))\in\R^{n+1}:q\in
C^{(\ell)}\}$. Then, for a sufficiently generic lifting function
$\omega$, the set of cells $C$ of $\mathcal{A}$ satisfying the
conditions:
\renewcommand{\theenumi}{$\roman{enumi}$}
\begin{enumerate}
\item  $\dim(conv(\widehat C^{(1)}+\cdots+\widehat C^{(s)})) = n$,
\item $(\widehat C^{(1)},\ldots,\widehat C^{(s)})$ is a face of
$(\widehat{\mathcal{A}}^{(1)}, \ldots,
\widehat{\mathcal{A}}^{(s)})$ whose inner normal has positive last
coordinate,
\end{enumerate}
\renewcommand{\theenumi}{\arabic{enumi}}
is a fine--mixed subdivision of $\mathcal{A}$ (see \cite[Section
2] {HuSt95}). More precisely, we have the following result (cf.
\cite[Section 2]{HuSt95}):
\begin{lemma}
The lifting process associated to a lifting function $\omega$
yields a fine-mixed subdivision of $\mathcal{A}$ if the following
condition holds: for every $r_1,\dots , r_s \in \Z_{\ge 0}$ with
$\sum_{\ell = 1}^s r_\ell > n$ and every cell $(C^{(1)}, \dots,
C^{(s)})$ with $C^{(\ell)} := \{ q_{\ell,0}, \dots, q_{\ell,
r_\ell}\} \subset \mathcal{A}^{(\ell)}$ $(1\le \ell \le s)$, if
$$
V(C):= \left(
\begin{array}{cccc}
q_{1,1} - q_{1,0}  \\
\vdots \\
q_{1,r_1} - q_{1,0}  \\
\cdots   \\
\cdots  \\
q_{s,1} - q_{s,0}  \\
\vdots \\
q_{s,r_s} - q_{s,0}
\end{array}\right)
\ \hbox{and} \ V(\widehat C) := \left(
\begin{array}{ccccccccc}
q_{1,1} - q_{1,0} & \omega_1(q_{1,1}) - \omega_1(q_{1,0})\\
\vdots & \vdots \\
q_{1,r_1}-q_{1,0} & \omega_1(q_{1,r_1}) - \omega_1(q_{1,0})\\
\cdots  & \cdots  \\
\cdots  & \cdots  \\
q_{s,1} - q_{s,0} & \omega_s(q_{s,1})-\omega_s(q_{s,0}) \\
\vdots & \vdots \\
q_{s,r_s}-q_{s,0} & \omega_s(q_{s,r_s})-\omega_s(q_{s,0})
\end{array}\right),
$$
then
$${\rm rank}(V(C)) = n \ \Longrightarrow \
{\rm rank} (V(\widehat C)) = n+1. $$\end{lemma}
\begin{proof} Notice that (1)--(3) are automatically satisfied by
the set of cells defined  by conditions $(i)$--$(ii)$. Assume that
the condition of the statement of the lemma is met and consider a
cell $C = (C^{(1)}, \dots, C^{(s)})$ of $\mathcal{A}$ satisfying
conditions $(i)$ and $(ii)$ above. Being $\widehat C$ a
\emph{lower} facet of $\mathcal{A}$, the identity
$\dim(conv(C^{(1)}+\cdots +C^{(s)})) = \dim( conv(\widehat
C^{(1)}+\cdots+\widehat C^{(s)}))$ must hold. Write $C^{(\ell)} =
\{ q_{\ell,0}, \dots, q_{\ell, r_\ell}\}$ for $1\le \ell \le s$.
Then we have that ${\rm rank}(V(C)) = \dim (< q_{\ell, j} -
q_{\ell, 0}: 1\le \ell \le r_\ell,\, 1\le j \le r_j>) = \dim(conv(
C^{(1)}+\cdots+C^{(s)})) = n$ and ${\rm rank}(V(\widehat C)) =
\dim(conv(\widehat C^{(1)}+\cdots+\widehat C^{(s)})) = n$. Now,
the condition stated on $\omega$ implies that $\sum_{\ell = 1}^s
r_\ell \le n$ and, taking into account that the $\sum_{\ell = 1}^s
r_\ell$ many vectors $q_{\ell, j} - q_{\ell,0}$ $(1\le \ell \le
s,\ 1\le j \le r_\ell)$ span a linear space of dimension $n$, we
conclude that the equality $\sum_{\ell = 1}^s r_\ell  =  n$ holds,
which shows that condition (5) in the definition of a fine--mixed
subdivision is met. Moreover, as $\sum_{\ell = 1}^s
\dim(conv(C^{(\ell)})) \ge \dim(conv(C^{(1)} +\cdots +C^{(s)} ))$
for an arbitrary cell $C$, we see that $\dim(conv(C^{(\ell)})) =
r_\ell$ holds for every $1\le \ell \le s$, which implies that
condition (4) is also valid. This finishes the proof of the lemma.
\end{proof}

Note that the condition ${\rm rank}(V(\widehat C)) = n+1$ can be
restated as the non--vanishing of the maximal minors of the matrix
$V(\widehat C)$. Since ${\rm rank}(V(C)) = n$, these maximal
minors are nonzero linear forms in the unknown values
$\omega_\ell(q_{\ell,j})$ of the lifting function. Thus, if
$\mathcal{N}_\ell = \# \mathcal{A}^{(\ell)}$ for every $1\le \ell
\le s$, a sufficiently generic lifting function can be obtained by
randomly choosing the values $\omega_\ell(q_{\ell,j})$ of $\omega$
at the points of $\mathcal{A}^{(\ell)}$ from the set
$\{1,2,\ldots,\rho 2^{\mathcal{N}_1+\cdots+\mathcal{N}_s} \}$,
with probability of success at least $1 - 1/\rho$ for $\rho \in
\mathbb{N}$.

In the sequel, we shall assume that a sufficiently generic lifting
function and the induced fine--mixed subdivision of $\mathcal{A}$
are given.
%
%
\subsection{Complexity model and complexity estimates}
%
In this section we describe our computational model and briefly
mention efficient algorithms for some basic specific algebraic
tasks.
%
%
\subsubsection{Complexity model}
Algorithms in computational algebraic geometry are usually
described using the standard dense (or sparse) complexity model,
i.e. encoding multivariate polynomials by  means  of  the vector
of all (or of all nonzero) coefficients. Taking into account that
a generic $n$--variate polynomial of degree~$d$
has~$\binom{d+n}{n}=O(d^n)$ nonzero coefficients, we see that the
dense   representation of multivariate polynomials requires an
exponential size, and their manipulation usually requires an
exponential number of arithmetic operations with respect to the
parameters~$d$ and~$n$. In order to avoid this exponential
behavior, we are going to use an alternative encoding of input,
output and intermediate results of our computations by means of
\slps\ (cf. \cite{Heintz89}, \cite{Strassen90}, \cite{Pardo95},
\cite{BuClSh97}). A {\sf \slp}~$\beta$ in $\Q(X):={\Q}(\xon)$ is a
finite sequence of rational functions ${(f_1\klk f_k)}\in
{\Q}(X)^k$ such that for~${1\le i\le k}$, the function $f_i$ is an
element of the set~$\{\xon\}$, or an element of~$\Q$ (a {\em
parameter}), or there exist~${1\le i_1,i_2<i}$ such that
$f_i=f_{i_1}\circ_i\,f_{i_2}$ holds, where $\circ_i$ is one of the
arithmetic operations~${+,-,\times ,{\div}}$. The \slp~$\beta$ is
called {\sf division--free} if~$\circ_i$ is different
from~${\div}$ for~${1\le i\le k}$. A natural measure of the
complexity of $\beta$ is its {\sf time} or {\sf length} (cf.
\cite{Borodin93}, \cite{Savage98}), which is the total number of
arithmetic operations performed during the evaluation process
defined by $\beta$. We say that the \slp~$\beta$ {\sf computes} or
{\sf represents} a subset~$S$ of ${\Q}(X)$ if $S\subset \{f_1\klk
f_k\}$ holds.

Our model of computation is based on the concept of \slps.
However, a model of computation consisting {\em only} of \slps\ is
not expressive enough for our purposes. Therefore we allow our
model to include decisions and selections (subject to previous
decisions). For this reason we shall also consider {\sf
computation trees}, which are \slps\ with {\sf branchings}. Time
of the evaluation of a given computation tree is defined similarly
to the case of \slps\ (see e.g. \cite{Gathen86}, \cite{BuClSh97}
for more details on the notion of computation trees). \spar
%
%
\subsubsection{Probabilistic identity testing}
A difficult point in the manipulation of multivariate polynomials
given by \slps\ is the so--called {\sf identity testing problem}:
given two elements $f$ and $g$ of $\C[X]:={\C}[\xon]$, decide
whether $f$ and $g$ represent the same polynomial function on
$\C^n$. Indeed, all known deterministic algorithms solving this
problem have complexity at least~$\max\{\deg f,\deg
g\}^{\Omega(1)}$. In this article we are going to use {\em
probabilistic} algorithms to solve the identity testing problem,
based on the following result:
\begin{theorem}[\cite{LiNi83}, \cite{Schmidt76}]
\label{th: Zippel-Schwartz} Let $f$ be a nonzero polynomial of
$\C[X]$ of degree at most~$d$ and let $\mathcal{S}$ be a finite
subset of $\C$. Then the number of zeros of $f$ in $\mathcal{S}^n$
is at most $d(\#\mathcal{S})^{n-1}$.
\end{theorem}
For the  analysis of our algorithms, we  shall interpret the
statement of Theorem~\ref{th: Zippel-Schwartz} in terms  of
probabilities. More precisely, given a fix nonzero polynomial $f$
in $\C[\xon]$ of degree at most~$d$, we conclude from
Theorem~\ref{th: Zippel-Schwartz} that the probability of choosing
randomly a point $a\in\mathcal{S}^n$ such that~${f(a)=0}$ holds is
bounded from above by $d/ \#\mathcal{S}$ (assuming a uniform
distribution of probability on the elements of $\mathcal{S}^n$).
\spar
%
%
\subsubsection{Basic complexity estimates}
In order to estimate the complexity of our procedures we shall
frequently use the notation ${\sf M}(m):=m\log^2m\log\log m$. Here
and in the sequel $\log$ will denote logarithm in base 2. Let $R$
be a commutative ring of characteristic zero with unity. We recall
that the number of arithmetic operations in $R$ necessary to
compute the multiplication or division with remainder of two
univariate polynomials in $R[T]$ of degree at most $m$ is
$O\big({\sf M}(m)/\log(m)\big)$ (cf. \cite{GaGe99},
\cite{BiPa94}). Multipoint evaluation and interpolation of
univariate polynomials of $R[T]$ of degree $m$ at {\em invertible}
points $a_1\klk a_m\in R$ can be performed with $O\big({\sf
M}(m)\big)$ arithmetic operations in $R$ (see e.g.
\cite{BoLeSc03}).

If $R=k$ is a field, then we shall use algorithms based on the
Extended Euclidean Algorithm (EEA for short) in order to compute
the gcd or resultant of two univariate polynomials in $k[T]$ of
degree at most $m$ with $O\big({\sf M}(m)\big)$ arithmetic
operations in $k$ (cf. \cite{GaGe99}, \cite{BiPa94}). We use
Pad\'e approximation in order to compute the dense representation
of the numerator and denominator of a rational function $f=p/q\in
k(T)$ with $\max\{\deg p,\deg q\}\le m$ from its Taylor series
expansion up to order $2m$. This also requires $O({\sf M}(m))$
arithmetic operations in $k$ (\cite{GaGe99}, \cite{BiPa94}).

For brevity, we will denote by $\Omega$ the exponent that appears
in the complexity estimate $O(n^\Omega)$ for the multiplication of
two $(n\times n)$--matrices with coefficients in $\Q$. We remark
that the (theoretical) bound $\Omega<2.376$ is typically
impractical and we prefer to take $\Omega:=\log 7\sim 2.81$ (cf.
\cite{BiPa94}).
%
%
\subsection{Geometric solutions}\label{subsect: geometric sol}
The notion of a geometric solution of an algebraic variety was
first introduced in the works of Kronecker and K{\"o}nig in the
last years of the XIXth century. Nowadays, geometric solutions are
widely used in computer algebra as a suitable representation of
algebraic varieties, especially in the zero--dimensional case.

Let $\overline{K}$ denote an algebraic closure of a field $K$ of
characteristic zero, let $\A^n(\overline{K})$ be the
$n$--dimensional space $\overline{K}^n$ endowed with its Zariski
topology, and let $V = \{ \xi^{(1)}, \dots,\xi^{(D)} \}$ be a
zero--dimensional subvariety of $\A^n(\overline{K})$ defined over
$K$. A {\sf geometric solution} of $V$ consists of
\begin{itemize}
\item a linear form $u = u_1 X_1 + \dots + u_n X_n\in K[X]$ which
separates the points of $V$, i.e. satisfying $u (\xi^{(i )})\ne u
(\xi^{(k)})$ if $i \ne k$,
\item the minimal polynomial  $m_u:= \prod_{1 \le i \le D} (Y -
u(\xi^{(i)}))\in K[Y]$ of $u$ in $V$ (where $Y$ is a new
variable),
\item polynomials $w_1, \dots,w_n \in K[Y]$ with $\deg w_j< D$ for
every $1\le j \le n$ satisfying
$$ V = \{(w_1(\eta) ,\dots,w_n (\eta) ) \in
\overline{K}^n\, / \, \eta \in \overline{K} ,\
 m_u(\eta) = 0 \}.$$
\end{itemize}

In the sequel, we shall be given a polynomial system
$f_1=\dots=f_n=0$ of $n$--variate polynomials of $\Q[X]$ defining
a zero--dimensional affine variety $V\subset\A^n:=\A^n(\C)$. We
shall consider the system $f_1=\dots=f_n=0$ (symbolically)
``solved'' if we obtain a geometric solution of $V$ as defined
above.

This notion of geometric solution can be extended to
equidimensional varieties of positive dimension. For our purposes,
it will be sufficient to consider the case of an algebraic curve
defined over $\Q$.

Suppose that we are given a curve $V\subset \A^{n+1}$ defined by
polynomials $f_1\klk f_n\in\Q[X,T]$. Assume that for each
irreducible component $C$ of $V$, the identity $I(C) \cap \Q[T] =
\{ 0 \}$ holds. Let $u$ be a nonzero linear form of $\Q[X]$ and
$\pi_u:V\to\A^2$ the morphism defined by $\pi_u(x,t):=(t,u(x))$.
Our assumptions on $V$ imply that the Zariski closure
$\overline{\pi_u(V)}$ of the image of $V$ under $\pi_u$ is a
hypersurface of $\A^{2}$ defined over $\Q$. Let $Y$ be a new
indeterminate. Then there exists a unique (up to scaling by
nonzero elements of $\Q$) polynomial $M_u\in\Q[T,Y]$ of minimal
degree defining $\overline{\pi_u(V)}$. Let $m_u\in\Q(T)[Y]$ denote
the (unique) monic multiple of $M_u$ with $\deg_Y (m_u)=\deg_Y
(M_u)$. We call $m_u$ the {\em minimal polynomial} of $u$ in $V$.
In these terms, a {\sf geometric solution} of the curve $V$
consists of
\begin{itemize}
\item a \emph{generic} linear form $u\in\Q[X]$,
\item the minimal polynomial $m_u\in\Q(T)[Y]$,
\item elements $v_1,\dots, v_n$ of $\Q(T)[Y]$ such that
$\frac{\partial m_u}{\partial Y}X_{i}\equiv v_{i} \mbox{ mod }
\Q(T)\otimes\Q[V]$ and $\deg_Y (v_{i})< \deg_Y (m_u )$ holds for
$1\le i\le n$.
\end{itemize}
We observe that $\deg_Y m_U$ equals the cardinality of the
zero-dimensional variety defined by $f_1,\dots, f_n$ over
$\A^n\big(\overline{\Q(T)}\big)$.

In the sequel, we shall deal with curves $V:=V(f_1\klk f_n)\subset
\A^{n+1}$ as above. The complexity of the algorithms for solving
the systems $f_1=\cdots=f_n=0$ defining such curves will be
expressed mainly by means of two discrete invariants: the {\em
degree} and the {\em height} of the projection $\pi:V\to\A^1$. The
degree of $\pi$ is defined as the degree $\deg m_u=\deg_YM_u$ of
the minimal polynomial of a generic linear form $u\in\Q[\xon]$ and
can be considered as a measure of the ``complexity'' of the curve
$V$. On the other hand, the height of $\pi$ is defined as
$\deg_TM_u$ and may be considered as a measure of the ``complexity
of the description'' of the curve $V$.

In the sparse setting, we can estimate $\deg_YM_u$ and $\deg_T
M_u$ in terms of combinatorial quantities (namely, mixed volumes)
associated to the polynomial system under consideration (see also
\cite{PhSo06}).
\begin{lemma}\label{lemma: estimate sparse height}
Let assumptions and notations be as above. For $1\le i\le n$, let
$Q_i\subset\R^n$ be the Newton polytope of $f_i$, considering
$f_i$ as an element of $\Q(T)[X]$. Let $\widehat Q_1,\dots,
\widehat Q_n\subset \R^{n+1}$ be the Newton polytopes of
$f_1,\dots, f_n$, considering $f_1\klk f_n$ as elements of
$\Q[X,T]$, and let $\Delta\subset \R^{n+1}$ be the standard
unitary simplex in the plane $\{T = 0\}$, i.e., the Newton
polytope of a generic linear form $u\in\Q[X]$. Assume that $0 \in
\widehat Q_i$ for every $1\le i \le n$. Then the following
estimates hold:
\begin{equation}\label{eq: estimate height polyhedral def}
\deg_Y M_u \le MV_n(Q_1,\dots, Q_n),\ 
\ \deg_T M_u\le MV_n(\Delta, \widehat Q_1,\dots, \widehat Q_n).
\end{equation}
Furthermore, if there exist $c_1\klk c_n\in \R_{\ge 0}$ such that
$\widehat Q_i \subset Q_i \times [0, c_i]$ for $1\le i\le n$, then
the following inequality holds:
\begin{equation}\label{eq: estimate height second def}
\deg_T M_u \le \sum_{i=1}^n c_i\, MV(\Delta, Q_1,\dots, Q_{i-1},
Q_{i+1},\dots,Q_n).\end{equation}
\end{lemma}

\begin{proof}
The upper bound $\deg_Y M_u \le MV_n(Q_1,\dots, Q_n)$ follows
straightforwardly from the BKK bound and the affine root count in
\cite{LiWa96}.


In order to obtain an upper bound for $\deg_TM_u$, we observe that
substituting a generic value $y\in \Q$ for $Y$ we have $\deg_T
M_u(T,\!Y)=\deg_T M_u(T, y)=\#\{t\in\C;M_u(t,y)=0\}$. Moreover, it
follows that $M_u(t,y)=0$ if and only if there exists a point
$x\in\A^n$ with $y = u(x)$ and $(x,t) \in V$. Thus, it suffices to
estimate the number of points $(x, t) \in \A^{n+1}$ satisfying
$u(x) - y = 0, f_1(x,t) = 0, \dots, f_n(x,t) = 0$. Being $u$ a
generic linear form, the system
\begin{equation}\label{eq: system for estimate height}
u(X) - y = 0, f_1(X,T) = 0, \dots, f_n(X,T)=0
\end{equation}
has finitely many common zeros in $\A^{n+1}$. Combining the BKK
bound with the affine root count of \cite{LiWa96} we see that
there are at most $MV(\Delta, \widehat Q_1,\dots, \widehat Q_n)$
solutions of (\ref{eq: system for estimate height}). We conclude
that $\deg_T M_u \le MV(\Delta, \widehat Q_1,\dots, \widehat Q_n)$
holds, which shows (\ref{eq: estimate height polyhedral def}).

In order to prove (\ref{eq: estimate height second def}), we make
use of basic properties of the mixed volume (see, for instance,
\cite[Ch. IV]{Ewald96}). Since $\widehat Q_i \subset Q_i \times
[0, c_i]$ holds for $1\le i\le n$, by the monotonicity of the
mixed volume we have
$$MV(\Delta, \widehat Q_1,\dots, \widehat Q_n) \le MV(\Delta,
Q_1\times [0,c_1],\dots, Q_n\times [0,c_n]).$$
Note that $Q_i\times [0,c_i] = S_{i,0} + S_{i,1}$, where $S_{i,0}
= Q_i \times \{ 0\}$ and $S_{i,1} =\{0\} \times [0,c_i] $ for
$i=1,\dots, n$. Hence, by multilinearity,
\begin{equation}\label{eq: MV}
MV(\Delta, Q_1\times [0,c_1],\dots, Q_n\times [0,c_n])=
\!\!\!\!\!\!\sum_{(j_1,\dots, j_n)\in
\{0,1\}^n}\!\!\!\!\!\!MV(\Delta, S_{1,j_1}, \dots, S_{n,j_n}).
\end{equation}

If the vector $(j_1,\dots, j_n)$ has at least two nonzero
coordinates, then two of the sets  $S_{1,j_1}, \dots, S_{n,j_n}$
are parallel line segments; therefore, $MV(\Delta, S_{1,j_1},
\dots, S_{n,j_n}) = 0$. On the other hand, if $j_i$ is the only
nonzero coordinate, the corresponding term in the sum of the
right--hand side of (\ref{eq: MV}) is
\begin{eqnarray*}
MV_{n+1}(\Delta, Q_1\times \{0\}, \dots, Q_{i-1}\times \{0\},
\{0\}\times [0,c_i], Q_{i+1}\times \{0\},\dots, Q_{n}\times \{0\})
\\[1mm] = c_i \,MV_n(\Delta, Q_1, \dots, Q_{i-1}, Q_{i+1}, \dots,
Q_{n}).\end{eqnarray*}
Finally, for $(j_1\dots, j_n) = (0,\dots, 0)$ we have
$MV_{n+1}(\Delta, Q_1\times \{0\},\dots, Q_n\times \{0\}) = 0$
since all the polytopes are included in an $n$-dimensional
subspace.

We conclude that the right-hand side of (\ref{eq: MV}) equals the
right--hand side of (\ref{eq: estimate height second def}). This
finishes the proof of the lemma.
\end{proof}

{From} the algorithmic point of view, the crucial step towards the
computation of a geometric solution of the variety $V$ consists in
the computation of the minimal polynomial $m_u$ of a generic
linear form $u$ which separates the points of $V$. In the
remaining part of this section we shall show how we can derive an
algorithm for computing the entire geometric solution of a
zero--dimensional variety $V$ defined over $\Q$ from a procedure
for computing the minimal polynomial of a generic linear form $u$
(cf. \cite{AlBeRoWo96}, \cite{GiLeSa01}, \cite{Schost03}).

Let $\Lambda:=(\Lambda_1\klk\Lambda_n)$ be a vector of new
indeterminates and let $K := \Q(\Lambda)$. Denote by $I_K$ the
ideal in $K[X_1,\dots, X_n]$ which is the extension of the ideal
$I:=I(V) \subset \Q[X_1,\dots, X_n]$ of the zero-dimensional
variety $V$, and denote by $B:=K[X_1,\dots, X_n]/I_K$ the
corresponding zero--dimensional quotient algebra. Write $V = \{
\xi^{(1)}, \dots, \xi^{(D)}\}$.

Set $U:=\Lambda_1X_1\plp \Lambda_nX_n\in K[X_1,\dots, X_n]$ and
let $m_U(\Lambda,Y) = \prod_{j=1}^D\big(Y - U(\xi^{(j)})\big)
\in\Q[\Lambda,Y]$ be the minimal polynomial of the linear form $U$
in the extension $K\hookrightarrow B$. Note that $\deg m_U = D$
holds, and that ${\partial m_U}/{\partial Y}$ is not a zero
divisor in $\Q[\A^n \times V]$. Furthermore, $m_U(\Lambda,{U}) \in
I(\A^n\times V) \subset \Q[\Lambda, X_1,\dots, X_n]$ holds. Since
$I(\A^n\times V)$ is generated by polynomials in $\Q[X_1,\dots,
X_n]$, taking the partial derivative of $m_U(\Lambda,{U})$ with
respect to the variable $\Lambda_k$ for $1\le k\le n$, we conclude
that
\begin{equation}
\label{eq: parametrizations} \frac{\partial m_U}{\partial Y}
(\Lambda,{U})\, X_k+\frac{\partial m_U}{\partial \Lambda_k}
(\Lambda,{U})\in I(\A^n\times V).
\end{equation}
Observe that the degree estimate $\deg_Y(\partial m_U/\partial
\Lambda_k)\le D-1$ holds.

Assume that a linear form $u = u_1 X_1 + \cdots + u_n X_n \in
\Q[X_1,\dots, X_n]$ which separates the points of $V$ is given.
Substituting $u_k$ for $\Lambda_k$ in the polynomial $m_U(\Lambda,
Y)$ we obtain the minimal polynomial $m_u(Y)$ of $u$. Furthermore,
making the same substitution in the polynomials $({\partial
m_U}/{\partial Y}) (\Lambda,Y)\, X_k+({\partial m_U}/{\partial
\Lambda_k}) (\Lambda,Y)$ of (\ref{eq: parametrizations}) for $1\le
k \le n$ and reducing modulo $m_{u}(Y)$, we obtain polynomials
$({\partial m_{u}}/{\partial Y})(Y)\, X_k-v_k(Y)\in I(V)
\quad(1\le k\le n)$. In particular, we have that the identities
\begin{equation}\label{eq: parametrizations curve}
\frac{\partial m_{u}}{\partial Y}(u)\, X_k = v_k(u)\ (1\le k \le
n)\end{equation}
hold in $\Q[V]$. Observe that the minimal polynomial $m_{u}(Y)$ is
square--free, since the linear form $u$ separates the points of
$V$. Therefore, $m_{u}(Y)$ and $\partial m_{u}/\partial Y(Y)$ are
relatively prime. Thus, multiplying modulo $m_u(Y)$ the
polynomials $v_k(Y)$ by the inverse of $(\partial m_{u}/\partial
Y)(Y)$ modulo $m_{u}(Y)$ we obtain polynomials $w_k(Y):=(\partial
m_{u}/\partial Y)^{-1}v_k(Y)$ $(1\le k\le n)$ of degree at most
$D-1$ such that
\begin{equation}\label{eq: parametrizations dimension zero}
X_k=w_k(u)\end{equation}
holds in $\Q[V]$ for $1\le k\le n$. The polynomials $m_u,w_1\klk
w_n\in\Q[Y]$ form a geometric solution of $V$.

Now, suppose that we are given an algorithm $\Psi$ over
$\Q(\Lambda)$ for computing the minimal polynomial of the linear
form $U=\Lambda_1 X_1 +\cdots +\Lambda_n X_n$. Suppose further
that we are given a separating linear form
$u:=u_1X_1+\dots+u_nX_n\in\Q[\xon]$ such that the vector $(u_1\klk
u_n)$ does not annihilate any denominator in $\Q[\Lambda]$ of any
intermediate result of the algorithm $\Psi$. In order to compute
the polynomials $v_1\klk v_n$ of (\ref{eq: parametrizations
curve}), we observe that the Taylor expansion of $m_U(\Lambda,Y)$
in powers of $\Lambda-u:=(\Lambda_1-u_1\klk\Lambda_n-u_n)$ has the
following expression:
$$m_U(\Lambda,Y)=m_{u}(Y)+\sum_{k=1}^n\Big( \frac{\partial
m_u}{\partial Y}(Y)X_k-v_k(Y)\Big)(\Lambda_k-u_k)\quad
\mathrm{mod}(\Lambda-u)^2.$$
We shall compute this first--order Taylor expansion by computing
the first--order Taylor expansion of each intermediate result in
the algorithm $\Psi$ . In this way, each arithmetic operation in
$\Q(\Lambda)$ arising in the algorithm $\Psi$ becomes an
arithmetic operation between two polynomials of $\Q[\Lambda]$ of
degree at most 1, and is truncated up to order $(\Lambda-u)^2$.
Since the first--order Taylor expansion of an addition,
multiplication or division of two polynomials of $\Q[\Lambda]$ of
degree at most 1 requires $O(n)$ arithmetic operations in $\Q$, we
have that the whole step requires $O(n{\sf T})$ arithmetic
operations in $\Q$, where ${\sf T}$ is the number of arithmetic
operations in $\Q(\Lambda)$ performed by the algorithm $\Psi$.

Finally, the computation of the polynomials $w_1,\dots, w_n$ of
(\ref{eq: parametrizations dimension zero}) requires the inversion
of $\partial m_{u}/\partial Y$ modulo $m_u(Y)$ and the modular
multiplication $w_k(Y):=(\partial m_{u}/\partial Y)^{-1}v_k(Y)$
for $1\le k\le n$. These steps can be executed with additional
$O\big(n{\sf M}(D)\big)$ arithmetic operations in $\Q$.
Summarizing, we have the following result:
\begin{lemma}\label{lemma: min/resgeom}
Suppose that we are given:
\begin{enumerate}
    \item an algorithm $\Psi$ in $\Q(\Lambda)$ which computes the minimal
polynomial $m_U\in\Q[\Lambda,Y]$ of $U:=\Lambda
X_1+\dots+\Lambda_nX_n$ with ${\sf T}$ arithmetic operations in
$\Q(\Lambda)$,
    \item a separating linear form $u:=u_1X_1+\dots+u_nX_n\in\Q[\xon]$
    such that the vector $(u_1\klk u_n)$ does not annihilate any
denominator in $\Q[\Lambda]$ of any intermediate result of the
algorithm $\Psi$.
\end{enumerate}
Then a geometric solution of the variety $V$ can be
(deterministically) computed with $O\big(n({\sf T}+{\sf
M}(D))\big)$ arithmetic operations in $\Q$.
\end{lemma}
%
%
\section{Statement of the problem and outline of the algorithm}
\label{sect: outline algoritmo}
Let $\Delta_1,\ldots,\Delta_n$ be fixed finite subsets of $\Z_{\ge
0}^n$ with $0\in\Delta_i$ for $1\le i\le n$ and let
$D:=MV(Q_1,\ldots,Q_n)$ denote the mixed volume of the polytopes
$Q_1:=conv(\Delta_1),\ldots,Q_n:=conv(\Delta_n)$. Assume that
$D>0$ holds or, equivalently, that  $\dim\left(\sum_{i\in I}
Q_i\right) \ge | I |$ for every non-empty subset $I\subset \{1,
\dots, n\}$ (see, for instance, \cite[Chapter IV, Proposition
2.3]{Oka97}).

Let $f_1,\ldots,f_n\in\Q[X]$ be polynomials defining a sparse
system with respect to $\Delta_1,\ldots,\Delta_n$ and let
$d_1,\ldots,d_n$ be their total degrees. Let $d:=\{d_1\klk d_n\}$.
Suppose that $f_1,\ldots,f_n$ define a zero--dimensional variety
$V$ in $\A^n$. As in the previous section, we group equal supports
into $s\le n$ distinct supports
$\mathcal{A}^{(1)},\ldots,\mathcal{A}^{(s)}$. Write
$\mathcal{A}:=(\mathcal{A}^{(1)},\ldots,\mathcal{A}^{(s)})$ and
denote by $k_\ell$ the number of polynomials $f_i$ with support
$\mathcal{A}^{(\ell)}$ for $1\le \ell\le s$.

{From} now on we assume that we are given a sufficiently generic
lifting function $\omega:=(\omega_1\klk\omega_s)$ and the
fine--mixed subdivision of $\mathcal{A}$ induced by $\omega$. We
assume further that the function $\omega_\ell:
\mathcal{A}^{(\ell)} \to \Z$ takes only nonnegative values and
$\omega_\ell(0\klk 0)=0$ for every $1\le \ell \le s$. The lifting
function $\omega$ and the corresponding fine--mixed subdivision of
$\mathcal{A}$ can be used in order to define an appropriate
deformation of the input system, the so-called \emph{polyhedral
deformation} introduced by Huber and Sturmfels in \cite{HuSt95}.
Our purpose here is to use this polyhedral deformation to derive a
symbolic probabilistic algorithm which computes a geometric
solution of the input system $f_1 = 0, \dots, f_n = 0$.

Since the polyhedral deformation requires that the coefficients of
the input polynomials satisfy certain generic conditions, we
introduce  some auxiliary {\em generic} polynomials $g_1\klk g_n$
with the same supports $\Delta_1\klk\Delta_n$ and consider the
perturbed polynomial system defined by $h_i:=f_i +g_i$ for $1\le i
\le n$. The genericity conditions underlying the choice of
$g_1\klk g_n$ and $h_1\klk h_n$ are discussed in Section
\ref{subsect: polyhedral deform}. We observe that if the
coefficients of the input polynomials $f_1,\ldots,f_n$ satisfy
these conditions then our method can be directly applied to
$f_1\klk f_n$.

Otherwise, we first solve the system $h_1 = 0, \dots, h_n = 0$ and
then recover the solutions to the input system $f_1 = 0, \dots,
f_n = 0$ by considering the homotopy $f_1 + (1-T) g_1=\dots=f_n +
(1-T) g_n=0$.
%
%
\subsection{The polyhedral deformation}
\label{subsect: polyhedral deform}
This section is devoted to introducing the polyhedral deformation
of Huber and Sturmfels.

Let $h_i:=\sum_{q\in\Delta_i}c_{i,q}X^q$ for $1\le i\le n$ be
polynomials in $\Q[X]$ and let $V_1$ denote the set of their
common zeros in $\A^n$. For $i=1, \dots, n$, let $\ell_i$ be the
(unique) integer with $\Delta_i = \mathcal{A}^{(\ell_i)}$, and let
$\widetilde\omega_i := \omega_{\ell_i}$ be the lifting function
associated to the support $\Delta_i$. In order to simplify
notations, the $n$--tuple $\widetilde\omega:=
(\widetilde\omega_1,\dots, \widetilde\omega_n)$ will be denoted
simply by $\omega = (\omega_1, \dots, \omega_n)$. As before, we
denote by
$\widehat{C}^{(\ell)}:=\{(q,\omega_\ell(q))\in\R^{n+1}:q\in
C^{(\ell)}\}$ the graph of any subset $C^{(\ell)}$ of
$\mathcal{A}^{(\ell)}$ for $1\le\ell\le s$, and extend this
notation correspondingly. For a new indeterminate $T$, we deform
the polynomials $h_1\klk h_n$ into polynomials
$\widehat{h}_1\klk\widehat{h}_n\in\Q[X,T]$ defined in the
following way:
\begin{equation}\label{eq: def h hat}\widehat{h}_i(X,T):=
\sum_{q\in\Delta_i}c_{i,q} X^qT^{\omega_i(q)} \quad (1\le i\le n).
\end{equation}
Let $I$ denote the ideal of $\Q[X,T]$ generated by
$\widehat{h}_1\klk\widehat{h}_n$ and let $J$ denote the Jacobian
determinant of $\widehat{h}_1\klk \widehat{h}_n$ with respect to
the variables $\xon$. We set
\begin{equation}\label{eq: def V hat}
\widehat V:=V(I:J^\infty)\subset\A^{n+1}.
\end{equation}
We shall show that, under a generic choice of the coefficients of
$h_1\klk h_n$, the system defined by the polynomials in (\ref{eq:
def h hat}) constitutes a deformation of the input system
$h_1=0\klk h_n=0$, in the sense that the morphism $\pi:\widehat
V\to\A^1$ defined by $\pi(x,t):=t$ is a dominant map with
$\pi^{-1}(1)= V_1\times \{1\}$. We shall further exhibit degree
estimates on the genericity condition underlying such choice of
coefficients. These estimates will allow us to obtain suitable
polynomials $h_1\klk h_n$ by randomly choosing their coefficients
in an appropriate finite subset of $\Z$.

According to \cite[Section 3]{HuSt95}, the solutions over an
algebraic closure $\overline{\Q(T)}$ of $\Q(T)$ to the system
defined by the polynomials (\ref{eq: def h hat}) are algebraic
functions of the parameter $T$ which can be represented as Puiseux
series of the form
\begin{equation}\label{eq: def series puiseux solucion de h gamma}
x(T) := (x_{10} T^{\gamma_1/\gamma_{n+1}} + \text{h-o t.}, \dots,
x_{n0} T^{\gamma_n/\gamma_{n+1}} + \text{h-o t.}),
\end{equation}
where $\gamma:=(\gamma_1,\dots, \gamma_n, \gamma_{n+1})\in
\Z^{n+1}$ is an inner normal with positive last coordinate
$\gamma_{n+1} >0$ of a (lower) facet $\widehat C = (\widehat
C^{(1)}, \dots, \widehat C^{(s)})$ of $\widehat{\mathcal{A}}$ of
type $(k_1,\dots, k_s)$, and $x_0:=(x_{10}, \dots, x_{n0}) \in
(\C^*)^n$ is a solution to the polynomial system defined by
\renewcommand{\theequation}{$\arabic{section}.\arabic{equation}_\gamma$}
\begin{equation}\label{eq: def h gamma 0}h_{i, \gamma}^{(0)} :=
\sum_{q\in C^{(\ell_i)}} c_{i, q} \, X^q \qquad (1\le i \le n),
\end{equation}
\renewcommand{\theequation}{\arabic{section}.\arabic{equation}}
$\!\!$where, as defined before, $\ell_i$ is the integer with $1\le
\ell_i\le s$ and $\Delta_i = \mathcal{A}^{(\ell_i)}$. We shall
``lift'' each solution $x_0$ to this system to a solution of the
form (\ref{eq: def series puiseux solucion de h gamma}) to the
system defined by (\ref{eq: def h hat}). This means that, on input
$x_0$, we shall compute the Puiseux series expansion of the
corresponding solution (\ref{eq: def series puiseux solucion de h
gamma}) truncated up to a suitable order.

Let
\begin{equation}\label{eq: V0gamma}
V_{0, \gamma}:= \{ x \in (\C^*)^n : h_{1,\gamma}^{(0)}(x) = 0, \dots,
h_{n,\gamma}^{(0)}(x) = 0\}.
\end{equation}
A particular feature of the polynomials (\ref{eq: def h gamma 0})
which makes the associated equation system ``easy to solve'' is
that the vector of their supports is $(C^{(1)})^{k_1}\times \cdots
\times (C^{(s)})^{k_s}$, where $(C^{(1)}, \dots, C^{(s)})$ is a
cell of type $(k_1,\dots, k_s)$ in a fine--mixed subdivision of
$\mathcal{A}$. Therefore, for every $1\le \ell\le s$, the set
$C^{(\ell)}$ consists of $k_\ell+1$ points and hence, up to
monomial multiplication so that each polynomial has a non-zero
constant term, the (Laurent) polynomials in (\ref{eq: def h gamma
0}) are linear combinations of $n+1$ distinct monomials in $n$
variables.

Denote $\Gamma\subset \Z^{n+1}$ the set of all primitive integer
vectors of the form $\gamma:= (\gamma_1,\dots,\! \gamma_n,
\!\gamma_{n+1}) \in \Z^{n+1}$ with $\gamma_{n+1} >0$ for which
there is a cell $C=(C^{(1)},\ldots,C^{(s)})$ of type $(k_1,\dots,
k_s)$ of the subdivision of $\mathcal{A}$ induced by $\omega$ such
that $\widehat C$ has inner normal $\gamma$.

Fix a cell $C=(C^{(1)},\ldots,C^{(s)})$ of type $(k_1,\dots, k_s)$
of the subdivision of $\mathcal{A}$ induced by $\omega$ associated
with a primitive inner normal $\gamma\in \Gamma$ with positive
last coordinate. In order to lift the points of the variety
$V_{0,\gamma}$ of (\ref{eq: V0gamma}) to a solution of the system
defined by the polynomials in (\ref{eq: def h hat}), we will work
with a family of auxiliary polynomials $h_{1,\gamma}\klk
h_{n,\gamma}\in\Q[X,T]$ which we define as follows:
\begin{equation}\label{eq: def h gamma}h_{i,\gamma}(X,T):=
T^{-m_i}\widehat{h}_i(T^{\gamma_1}X_1\klk
T^{\gamma_n}X_n,T^{\gamma_{n+1}}) \quad (1\le i\le
n)\end{equation}
where $m_i\in\Z$ is the lowest power of $T$ appearing in
$\widehat{h}_i(T^{\gamma_1}X_1\klk
T^{\gamma_n}X_n,T^{\gamma_{n+1}})$ for every $1\le i \le n$. Note
that the polynomials obtained by substituting $T = 0$ into
$h_{1,\gamma}, \dots, h_{n, \gamma}$ are precisely those
introduced in (\ref{eq: def h gamma 0}).
%
%
\subsection{On the genericity of the initial system}
\label{subsect: genericity conditions}
Here we discuss the genericity conditions underlying the choice of
the polynomials $g_1\klk g_n$ that enable us to apply the
polyhedral deformation defined by the lifting form $\omega$ to the
system $h_1:= f_1+g_1=0, \dots, h_n := f_n +g_n=0$.

The first condition we require is that   the set of common zeros
of the perturbed polynomials $h_1 ,\dots, h_n$ is a
zero--dimensional variety with the maximum number of points for a
sparse system with the given structure. More precisely, we require
the following condition:
\begin{itemize}
\item[$(\mathsf{H1})$] The set $V_1:=\{x\in\A^n:h_1(x) = 0, \dots,
h_n(x) = 0\}$ is a zero-dimensional variety with
$D:=MV(Q_1,\dots,Q_n)$ distinct points.
\end{itemize}

In addition,  we need that the system (\ref{eq: def h gamma 0})
giving the initial points to our first deformation for every
$\gamma\in\Gamma$ has as many roots as possible, namely the mixed
volume of their support vectors.

{For} each cell $C = (C^{(1)}, \dots, C^{(s)}) $ of type
$(k_1,\dots, k_s)$ of the induced fine--mixed subdivision, set an
order on the $n+1$ points appearing in any of the sets
$C^{(\ell)}$, after a suitable translation so that $0 \in
C^{(\ell)}$ for every $1\le \ell \le s$. Assume that $0\in \Z^n$
is the last point according to this order. Denote $\gamma\in
\Z^{n+1}$ the primitive inner normal of $C$ with positive last
coordinate. Consider  the $n\times (n+1)$ matrix whose $i$th row
is the coefficient vector of $h_{i,\gamma}^{(0)}$ in the
prescribed monomial order and set $\mathcal{M}_{\gamma}\in
\Q^{n\times n}$ and $\mathcal{B}_{\gamma}\in \Q^{n\times 1}$ for
the submatrices consisting of the first $n$ columns (coefficients
of non-constant monomials) and the last column (constant
coefficients) respectively. Then, the coefficients of $g_1,\dots,
g_n$ are to be chosen so that the following condition holds:
\begin{itemize}
\item[$\mathsf{(H2)}$] For every $\gamma \in \Gamma$, the
$(n\times n)$--matrix $\mathcal{M}_{\gamma}$ is nonsingular and
all the entries of
$(\mathcal{M}_{\gamma})^{-1}\mathcal{B}_{\gamma}$ are nonzero.
\end{itemize}
\medskip

Our next results assert that the above conditions can be met with
good probability by randomly choosing the coefficients of
$g_1,\ldots,g_n$ in a certain set $\mathcal{S}\subset\Z$. We
observe that our estimate on the size of $\mathcal{S}$ is not
intended to be accurate, but to show that the growth of the size
of the integers involved in the subsequent computations is not
likely to create complexity problems.

Let $\{\Omega_{i,q}:  1\le i\le n,\ q\in\Delta_i\}$ be a set of
new indeterminates over $\Q$. For $1\le i \le n$, write
$\Omega_i:=(\Omega_{i,q}:{q\in\Delta_i})$ and let $H_i\in
\Q[\Omega_i, X]$ be the generic polynomial
\begin{equation}\label{eq: His}
H_i(\Omega_i,X):=\sum_{q\in\Delta_i}\Omega_{i,q}X^q
\end{equation}
with support $\Delta_i$ and $N_i:= \# \Delta_i$ coefficients. Let
$\Omega:=(\Omega_1,\ldots,\Omega_n)$ and let $N:= N_1 +\cdots
+N_n$ be the total number of indeterminate coefficients.

We start the analysis of the required generic conditions with the
following quantitative version of Bernstein's result on the
genericity of zero-dimensional sparse systems (see \cite[Theorem
B]{Bernstein75}, \cite[Theorem 6.1]{HuSt95}):

\begin{lemma}\label{lemma: sparse system is 0dim}
There exists a nonzero polynomial $P^{(0)}\in\Q[\Omega]$ with
$\deg P^{(0)} \le 3 n^{2n+1} d^{2n-1}$ such that for any $c
\in\Q^N$ with $P^{(0)}(c)\not=0$, the system $H_1(c_1,X)=0, \dots,
H_n(c_n,X)=0$ has $D$ solutions in $\A^n$, counting
multiplicities.
\end{lemma}
\begin{proof}
Due to \cite[Theorem 6.1]{HuSt95} combined with \cite{LiWa96}, the
system $H_1(c_1,X)=0, \dots$, $H_n(c_n,X)=0$ has $D$ solutions in
$\A^n$ counting multiplicities if and only if for every facet
inner normal $\mu\in \Z^n$ of $Q_1+\cdots +Q_n$, the sparse
resultant ${\rm Res}_{\Delta_1^\mu, \dots,\Delta_n^\mu}$ does not
vanish at $c:=(c_1,\dots, c_n)$. Here $\Delta_i^\mu$ denotes the
set of points of $\Delta_i$ where the linear functional induced by
$\mu$ attains its minimum for $1\le i \le n$.

Therefore, the polynomial $P^{(0)}:= \prod_{\mu}{\rm
Res}_{\Delta_1^\mu, \dots,\Delta_n^\mu} \in \Q[\Omega]$, where the
product ranges over all primitive inner normals $\mu \in \Z^n$ to
facets of $Q_1 +\cdots + Q_n$, satisfies the required condition.

In order to estimate the degree of $P^{(0)}$, we observe that for
every facet inner normal $\mu\in \Z^n$ the following upper bound
holds:
$$\deg({\rm Res}_{\Delta_1^\mu, \dots,\Delta_n^\mu}) \le
\sum_{i=1}^n MV(\Delta_1^\mu, \dots, \widehat{\Delta_i^\mu},
\dots, \Delta_n^\mu)\le n d^{n-1},$$ where $d:=\max\{d_1,\dots,
d_n\}$. On the other hand, it is not difficult to see that the
number of facets of an $n$-dimensional integer convex polytope
$P\subset \R^n$ which has an integer point in its interior is
bounded by $n!\, vol_{\R^n}(P)$. Now, taking $P:= (n+1) Q$, we
obtain an integer polytope with the same number of facets as $Q$
having an integer interior point. Then, the number of facets of
$Q$ is bounded by $n! \, vol_{\R^n}(P) = n!\, vol_{\R^n}((n+1)Q) =
(n+1)^n\, n!\, vol_{\R^n}(Q) \le (n+1)^n (nd)^n$, since $Q$ is
included in the $n$-dimensional simplex of size $nd$. This proves
the upper bound for the degree $P^{(0)}$ of the statement of the
lemma.
\end{proof}

The next lemma is concerned with the genericity of a smooth sparse
system.

\begin{lemma} \label{lemma: grado condicion H1} With the same
notations as in Lemma \ref{lemma: sparse system is 0dim} and
before, there exists a nonzero polynomial $P^{(1)}\in\Q[\Omega]$
of degree at most $4 n^{2n+1} d^{2n-1}$ such that for any
$c\in\Q^N$ with $P^{(1)}(c)\not=0$, the system $H_1(c_1,X)=0,
\dots,H_n(c_n,X)=0$ has exactly $D$ distinct solutions in
$\A^n$.\end{lemma}

\begin{proof}
Consider the incidence variety associated to
$(\Delta_1\klk\Delta_n)$--sparse systems, namely
$$W:=\{(x,c)\in(\C^*)^n\times(\A^{N_1}\times\cdots\times \A^{N_n}):\sum_{q\in
\Delta_i}\!\!c_{i,q}x^q=0\ \mathrm{for}\ 1\le i\le n\}.$$
As in \cite[Proposition 2.3]{PeSt93}, it follows that $W$ is a
$\Q$-irreducible variety. Let $\pi_\Omega: W \to \A^{N_1} \times
\cdots \times \A^{N_n}$ be the canonical projection, which is a
dominant map.

By \cite[Chapter V, Corollary (3.2.1)]{Oka97}, there is a nonempty
Zariski open set $\mathcal{U}(\Delta_1,\dots, \Delta_n) \subset
\A^{N_1}\times \dots \times \A^{N_n}$ of coefficients $c =
(c_1,\dots, c_n)$ for which the polynomials $H_1(c_1,X),
\dots,H_n(c_n,X)$ have supports $\Delta_1, \dots,\Delta_n$
respectively and the set of their common zeros in $(\C^*)^n$ is a
non--degenerate complete intersection variety. Then, the Jacobian
$J_H:=\det(\partial H_i/\partial X_j)_{1\le i,j\le n}$ does not
vanish at any point of $\pi_\Omega^{-1}(c)$ for every $c\in
\mathcal{U}(\Delta_1, \dots, \Delta_n)$.

Let $\Q(\Omega)\hookrightarrow \Q(W)$ be the finite field
extension induced by the dominant projection $\pi_\Omega$. By the
preceding paragraph we have that the rational function defined by
$J_H$ in $\Q(W)$ is nonzero. Therefore, its primitive minimal
polynomial $M_J\in\Q[\Omega,T]$ is well defined and satisfies the
degree estimates
$$\deg_\Omega M_J\le \deg W \cdot \deg J_H \le \prod_{i=1}^n
(d_i+1) \cdot \ \sum_{i=1}^n d_i\le 2^n d^{n+1} n$$
(see \cite{SaSo96}, \cite{Schost03}).

Let $P^{(1)}:=P^{(0)}M_J^{(0)}$, where $P^{(0)}$ is the polynomial
given by Lemma \ref{lemma: sparse system is 0dim} and $M_J^{(0)}$
denotes the constant term of the expansion of $M_J$ in powers of
$T$. We claim that $P^{(1)}$ satisfies the requirements of the
statement of the lemma. Indeed, let $c\in\Q^N$ satisfy
$P^{(1)}(c)\not=0$. Then $P^{(0)}(c)\not=0$ holds and so, Lemma
\ref{lemma: sparse system is 0dim} implies that
$H_1(c,X)=\dots=H_n(c,X)=0$ is a zero-dimensional system.
Furthermore, $M_J^{(0)}(c)$ is a nonzero multiple of the product
$\prod_{x\in\pi_\Omega^{-1}(c)}J_H(c,x)$. Thus, the non-vanishing
of $M_J^{(0)}(c)$ shows that all the points of
$\pi_\Omega^{-1}(c)$ are smooth and therefore, from e.g. \cite[IV,
Theorem 2.2]{Oka97}, it follows that $\pi_\Omega^{-1}(c)$ consists
of exactly $D$ simple points in $(\C^*)^n$. Moreover, combining
the assumption that $0\in\Delta_i$ for $1\le i\le n$ with
\cite[Theorem 2.4]{LiWa96}, we deduce that $\pi_\Omega^{-1}(c)$
consists of $D$ simple points in $\A^n$. The estimate $\deg
M_J^{(0)}\le \deg_\Omega M_J\le 2^n d^{n+1} n \le n^{2(n+1)}
d^{2n-1}$ implies the statement of the lemma.
\end{proof}

Finally, we exhibit a generic condition on the coefficients
$h_1,\dots, h_n$ which implies that assumption $\mathsf{(H2)}$
holds.

\begin{lemma}\label{lemma: grado condicion H2}
With the previous assumptions and notations, there exists a nonzero
polynomial $P^{(2)}\in\Q[\Omega]$ with $\deg P^{(2)}\le
n(n+1)\#\Gamma$ such that for every $c:=(c_1,\ldots,c_n)\in\Q^N$
with $P^{(2)}(c)\neq 0$, the polynomials $h_i:= H_i(c_i,X)$ $(1\le
i\le n)$ satisfy condition $\mathsf{(H2)}$.
\end{lemma}
\begin{proof}
Fix a primitive integer inner normal $\gamma\in \Gamma$ to a lower
facet of   $\widehat{\mathcal{A}}$. Let $\mathcal{M}_{\gamma}\in
\Q[\Omega]^{n\times n}$ and $\mathcal{B}_{\gamma}\in
\Q[\Omega]^{n\times 1}$ be the matrices constructed from the
generic polynomials $H_1,\dots, H_n \in \Q[\Omega][X]$ as
explained in the paragraph preceding condition $\mathsf{(H2)}$.
Let $D_{0, \gamma} \in \Q[\Omega]$ be the (non-zero) determinant
of $\mathcal{M}_\gamma$, and for every $1\le j \le n$, let $D_{j,
\gamma}$ be the determinant of the matrix obtained from
$\mathcal{M}_\gamma$ by replacing its $j$th column with
$\mathcal{B}_\gamma$. Set $P_\gamma:= \prod_{j=0}^n D_{j,
\gamma}$. Finally, take $P^{(2)} := \prod_{\gamma\in \Gamma}
P_\gamma$. By Cramer's rule, whenever $P^{(2)}(c) \ne 0$, we have
that the system $h_1,\dots, h_n$ with coefficient vector $c =
(c_1,\dots, c_n)$ meets condition $\mathsf{(H2)}$.

The degree estimate for $P^{(2)}$ follows from the fact that $\deg
P_\gamma \le n(n+1)$ holds for every $\gamma \in \Gamma$, since
each of the entries of the matrices whose determinants are
involved has degree $1$ in the variables $\Omega$.
\end{proof}

Now, we are ready to state a generic condition on the coefficients
of $h_1\klk h_n$ which implies that $\mathsf{(H1)}$ and
$\mathsf{(H2)}$ hold.

\begin{proposition}\label{prop: cond final preparation}
Under the previous assumptions and notations, there exists a
nonzero polynomial $P\in\Q[\Omega]$ with $\deg P \le 4 n^{2n+1}
d^{2n-1} +n(n+1)D$ such that for every $c\in\Q^N$ with $P(c)\neq
0$, the polynomials $h_i:= H_i(c_i,X)$ $(1\le i\le n)$ satisfy
conditions $\mathsf{(H1)}$ and $\mathsf{(H2)}$.
\end{proposition}
\begin{proof}
Set $P:= P^{(1)} P^{(2)}$, where $P^{(1)}$ is the polynomial of
the statement of Lemma \ref{lemma: grado condicion H1} and
$P^{(2)}$ is the one defined in the statement of Lemma \ref{lemma:
grado condicion H2}. The result follows from Lemmas \ref{lemma:
grado condicion H1} and \ref{lemma: grado condicion H2}, and the
upper bound $\# \Gamma \le D$ for the cardinality of the set of
the distinct inner normal vectors considered (one for each cell of
type $(k_1,\dots, k_s)$ in the given fine-mixed subdivision).
\end{proof}
%
%
\subsection{Outline of the algorithm}
Now we have all the tools necessary to give an outline of our
algorithm for the computation of a geometric solution of the
(sufficiently generic) sparse system $h_1=\dots=h_n=0$.

With notations as in the previous subsections, we assume that a
fine--mixed subdivision of $\mathcal{A}$ induced by a lifting
function $\omega$ is given. This means that we are given the set
$\Gamma$ of inner normals of the lower facets of the convex hull
of $\widehat{\mathcal{A}}$, together with the corresponding cells
of the convex hull of $\mathcal{A}$. In addition, we suppose that
our input polynomials $h_1, \dots, h_n \in \Q[X]$ satisfy
conditions $(\mathsf{H1})$ and $(\mathsf{H2})$ and denote by
$V_1\subset \A^n$   the affine variety defined by $h_1, \dots,
h_n$.

First, we choose a \emph{generic} linear form $u\in\Q[X]$ such
that:
\begin{itemize}
\item $u$ separates the points of the zero--dimensional varieties
$V_1$ and $V_{0,\gamma}$ for every $\gamma \in \Gamma$. This
condition is represented by the nonvanishing of a certain
nonconstant polynomial of degree at most $2D^2$.
\item An algorithm for the computation of the minimal polynomial
of $u$ in $V_{0,\gamma}$ to be described below can be extended to
a computation of a geometric solution of $V_{0,\gamma}$ according
to Lemma \ref{lemma: min/resgeom} for every $\gamma\in\Gamma$.
This condition is represented by the nonvanishing of a nonconstant
polynomial of degree at most $4D_\gamma^3$ for each
$\gamma\in\Gamma$.
\item An algorithm for the computation of the minimal polynomial
of $u$ in $\widehat{V}$ to be described below can be extended to a
computation of a geometric solution of $\widehat{V}$ according to
Lemma \ref{lemma: min/resgeom}. This application of Lemma
\ref{lemma: min/resgeom} requires that the coefficient vector of
the linear form $u$ does not annihilate a nonconstant polynomial
of degree at most $4D^4$.
\end{itemize}
Fix $\rho\ge 2$. From Theorem \ref{th: Zippel-Schwartz} it follows
that a linear form $u$ satisfying these conditions can be obtained
by randomly choosing its coefficients from the set $\{1\klk 6\rho
D^4\}$ with error probability at most $1/\rho$.

Next we compute the monic minimal polynomial $\widehat m_u\in
\Q(T)[Y]$ of the linear form $u$ in the curve $\widehat V$
introduced in (\ref{eq: def V hat}). For this purpose, we
approximate the Puiseux series expansions of the branches of
$\widehat V$ lying above 0  by means of a symbolic
(Newton--Hensel) ``lifting'' of the common zeros of the
zero--dimensional varieties $V_{0, \gamma}\subset \A^n$ defined by
the polynomials (\ref{eq: def h gamma 0}) for all
$\gamma\in\Gamma$ (see Section \ref{sect: sol V hat}).

 This in turn requires the computation of a geometric
solution of $V_{0, \gamma}$ for every $\gamma \in \Gamma$. By
means of a change of variables we put the system $h_{1,
\gamma}^{(0)}=\dots= h_{n,\gamma}^{(0)}=0$ defining the variety
$V_{0, \gamma}$ into a ``diagonal'' form (see Subsection
\ref{subsect: geo sol V0 gamma} below), which allows us to compute
the minimal polynomial $m_{u,\gamma}^{(0)}$ of $u$ in $V_{0,
\gamma}$. Since the linear form $u$ satisfies condition 2 of the
statement of Lemma \ref{lemma: min/resgeom}, from this procedure
we derive an algorithm computing a geometric solution of $V_{0,
\gamma}$ according to Lemma \ref{lemma: min/resgeom}.

Then we ``lift'' this geometric solution to a suitable
(non--archimedean) approximation $\widetilde m_{\gamma}$ of a
factor $m_{\gamma}$ (over $\overline{\Q(T)}$) of the desired
minimal polynomial $\widehat m_u$ of $u$. In the next step
 we obtain the minimal polynomial $\widehat m_u =
\prod_{\gamma\in \Gamma} m_{\gamma}$ from the approximate factors
$\widetilde m_{\gamma}$, namely, we compute the dense
representation of the coefficients (in $\Q(T)$) of
$\widehat{m}_u$, using Pad\'e approximation (see Subsection
\ref{subsect: comput geo sol V hat} below). Finally, we apply the
proof of Lemma \ref{lemma: min/resgeom} to derive an algorithm for
computing a geometric solution of the variety $\widehat{V}$.

In the last step we compute a geometric solution of the variety
$V_1$ by substituting 1 for $T$ in the polynomials that form the
geometric solution of $\widehat V$.

The whole algorithm for solving the system $h_1=\dots=h_n=0$ may
be briefly sketched as follows:
\begin{algorithm}
${}^{}$

\begin{itemize}
\item Choose the coefficients of a linear form $u\in\Q[X]$ at
random from the set $\{1\klk 6\rho D^4\}$.
\item For each $\gamma \in \Gamma:$
\begin{itemize}
\item Find a geometric solution of the variety $V_{0, \gamma}$
defined in {\rm (\ref{eq: V0gamma})}.
\item Obtain a straight-line program for the polynomials
$h_{1,\gamma}, \dots, h_{n,\gamma}$ defined in {\rm (\ref{eq: def
h gamma})} from the coefficients of $h_1,\dots, h_n$ and the
entries of $\gamma\in \Z^{n+1}$.
\item ``Lift'' the computed geometric solution of $V_{0, \gamma}$
to an approximation $\widetilde m_{\gamma}$ of the factor
$m_{\gamma}$ of $\widehat m_u$ by means of a symbolic
Newton--Hensel procedure.
\end{itemize}
\item Obtain a geometric solution of the curve $\widehat V:$
\begin{itemize}
    \item Compute the approximation $\widetilde{m}_u:=
    \prod_{\gamma\in\Gamma}\widetilde m_{\gamma}$ of
    $\widehat m_u$.
    \item Compute the dense representation of $\widehat m_u$
    from $\widetilde{m}_u$ using Pad\'e approximation.
    \item Find a geometric solution of $\widehat{V}$ applying
    the proof of Lemma \ref{lemma: min/resgeom}.
\end{itemize}
\item Substitute 1 for $T$ in the polynomials which form the
geometric solution of $\widehat V$ computed in the previous step
to obtain a geometric solution of the variety $V_1$.
\end{itemize}
\end{algorithm}
%
%
\section{Solution of the variety $\widehat V$}
\label{sect: sol V hat}
%
%
\subsection{Geometric solutions of the starting
varieties}\label{subsect: geo sol V0 gamma}
In this subsection we exhibit an algorithm that computes, for a
given inner normal $\gamma \in \Gamma$, a geometric solution of
the variety $V_{0, \gamma}\subset (\C^*)^n$ defined by the
polynomials $h_{i, \gamma}^{(0)}$ $(1\le i \le n)$ for polynomials
$h_1,\dots, h_n$ satisfying assumptions $\mathsf{(H1)}$ and
$\mathsf{(H2)}$. This algorithm is based on the procedure
presented in \cite{HuSt95}.

Fix a cell $C=(C^{(1)}\klk C^{(s)})$ of type $(k_1,\dots, k_s)$ of
the given fine--mixed subdivision of $\mathcal{A}$ and let $\gamma
\in \Gamma$ be its associated inner normal. For $1\le \ell\le s$,
we denote by $h_1^{(\ell)},\ldots,h_{k_\ell}^{(\ell)}$ the
polynomials in the set $\{h_{1,\gamma}^{(0)},\ldots,
h_{n,\gamma}^{(0)}\}$ that are supported in $C^{(\ell)}$. In the
sequel, whenever there is no risk of confusion we will not write
the subscript $\gamma$ indicating which cell we are considering.

Our hypotheses imply that $h_1^{(\ell)},\ldots,
h_{k_\ell}^{(\ell)}$ are $\Q$--linear combinations of precisely
$k_\ell+1$ monomials in $\Q[X]$ and, up to a multiplication by a
monomial, we may assume one of them to be the constant term.
Denote these monomials by
$X^{\alpha_{\ell,0}},\ldots,X^{\alpha_{\ell,k_\ell}}$, with
$\alpha_{\ell,0}:=0\in\Z^n$. Let
$\widetilde{\mathcal{M}}^{(\ell)}$ be the matrix of
$\Q^{k_\ell\times(k_\ell+1)}$ for which the following equality
holds in $\Q[X, X^{-1}]^{k_\ell}$:
\begin{equation}\label{eq: matrixC}
\widetilde{\mathcal{M}}^{(\ell)}\begin{pmatrix} X^{\alpha_{\ell,k_\ell}}\\
\vdots
\\X^{\alpha_{\ell,0}}\end{pmatrix} =
\begin{pmatrix}h_1^{(\ell)}\\ \vdots
\\h_{k_\ell}^{(\ell)}\end{pmatrix},
\end{equation}
and let $\mathcal{M}^{(\ell)}$ denote the square $(k_\ell\times
k_\ell)$--matrix obtained by deleting the last column from
$\widetilde{\mathcal{M}}^{(\ell)}$. Set
$$
\mathcal{M}:=\begin{pmatrix} \mathcal{M}^{(1)} & 0     & \cdots
& 0 \\ 0 & \mathcal{M}^{(2)} & \cdots & 0\\ 0     &  0     & \ddots & 0\\
0     & 0 & \cdots & \mathcal{M}^{(s)}
\end{pmatrix},$$
where $0$ here represents different block matrices with all its
entries equal to $0\in\Q$. Then $\mathcal{M}$ is the matrix
defined by the coefficients of the nonconstant terms of the
(Laurent) polynomials $h_{1,\gamma}^{(0)}, \ldots,
h_{n,\gamma}^{(0)}$, up to a translation.

Due to condition $\mathsf{(H2)}$ we have that the matrix
$\mathcal{M}$ is invertible, which in turn implies that the square
matrices $\mathcal{M}^{(\ell)}$ are invertible for $1\le \ell\le
s$. Following \cite{HuSt95}, we  apply Gaussian elimination to the
matrix $\widetilde{\mathcal{M}}^{(\ell)} $ for $1\le \ell\le s$
and obtain a set of $k_\ell+1$ binomials
$$
\begin{pmatrix}
1 & 0  & 0 & \ldots & -c_{\alpha_{\ell,k_\ell}}\\ 0 & 1 & 0 & \ldots
& -c_{\alpha_{\ell,k_\ell-1}}\\ \vdots & & & \ddots &  \\ 0 & 0
&\ldots & 1 & -c_{\alpha_{\ell,1}}
\end{pmatrix}
\begin{pmatrix}
X^{\alpha_{\ell,k_\ell}}\\ X^{\alpha_{\ell,k_\ell-1}} \\ \vdots
\\X^{\alpha_{\ell,1}}
\end{pmatrix} =
\begin{pmatrix}
X^{\alpha_{\ell,k_\ell}}-c_{\alpha_{\ell,k_\ell}} \\
X^{\alpha_{\ell,k_\ell-1}}-c_{\alpha_{\ell,k_\ell-1}}\\ \vdots
\\  X^{\alpha_{\ell,1}}-c_{\alpha_{\ell,1}} \end{pmatrix}$$
that generate the same linear subspace of $\Q[X, X^{-1}]$ as the
polynomials in (\ref{eq: matrixC}). Therefore, for $1\le\ell\le s$
the set of common zeros in $(\C^*)^n$ of the polynomials
$h_1^{(\ell)}\klk h_{k_\ell}^{(\ell)}$ is given by the system
$X^{\alpha_{\ell,k_\ell}}= c_{\alpha_{\ell,k_\ell}}, \ldots,
X^{\alpha_{\ell,1}}=c_{\alpha_{\ell,1}}$. Putting these $s$
systems together, we obtain a binomial system of the form
\begin{equation}\label{eq: sistema binomial}
X^{\alpha_1}= p_1, \ldots,  X^{\alpha_n}=p_n,
\end{equation}
with $\alpha_i\in\Z^n$ and $p_i\in\Q\setminus\{0\}$ $(1\le i \le
n)$, that defines the variety $V_{0,\gamma}$. Note that the second
part of condition $\mathsf{(H2)}$ ensures the non--vanishing of
the constants $p_i$ for $1\le i \le n$.

Now, let $\mathcal{E}$ denote the $(n\times n)$--matrix whose
columns are the exponent vectors $\alpha_1,\dots, \alpha_n$. Using
\cite[Proposition 8.10]{Storjohann00}, we obtain unimodular
matrices $K=(k_{i,j})_{1\le i,j\le n}$, $L=(l_{i,j})_{1\le i,j\le
n}$ of $\Z^{n\times n}$, and a diagonal matrix ${\rm
diag}(r_1,\ldots,r_n)\in\Z^{n\times n}$ which give the Smith
Normal Form for $\mathcal{E}$, i.e., matrices such that the
identity
\begin{equation}\label{eq: smith normal form matrix E}
K \cdot\mathcal{E}\cdot L={\rm diag}(r_1,\ldots,r_n)
\end{equation}
holds in $\Z^{n\times n}$. We observe that the upper bound
\begin{equation}\label{eq: estimate height Smith form}
\log\|K\|\le (4n+5)(\log n+\log\|\mathcal{E}\|)
\end{equation}
holds, where $\|A\|$ denotes the maximum of the absolute value of
the entries of a given matrix $A$ \cite[Proposition
8.10]{Storjohann00}.

Let $Z_1,\ldots,Z_n$ be new indeterminates, and write
$Z:=(Z_1,\ldots,Z_n)$. We introduce the change of coordinates
given by $X_i:=Z_1^{k_{1,i}}\cdots Z_n^{k_{n,i}}$ for $1\le i\le
n$. Making this change of coordinates in (\ref{eq: sistema
binomial}) we obtain the system
$$ Z^{K\alpha_1}= p_1, \ldots, Z^{K\alpha_n}=p_n, $$
which is equivalent to the ``diagonal'' system
$$
Z_j^{r_j} = \prod_{i=1}^n (Z^{K\alpha_i})^{l_{i,j}}= \prod_{i=1}^n
p_i^{l_{i,j}}  =:q_j\quad (1\le j\le n).
$$
Inverting some of the coefficients $q_j$ if necessary we may
assume without loss of generality that the integers $r_1,\dots,
r_n$ are positive.

We first describe an algorithm for computing a geometric solution
of the variety $V_{0,\gamma}\subset\A^n$ in the coordinate system
of $\A^n$ defined by $Z_1,\dots, Z_n$. This algorithm takes as
input the set of polynomials $Z_1^{r_1} - q_1,\dots, Z_n^{r_n} -
q_n \in \Q[Z_1,\dots, Z_n]$ defining $V_{0,\gamma}$ in the
coordinates $Z_1\klk Z_n$, and outputs a linear form
$\widetilde{u}\in\Q[Z_1\klk Z_n]$ which separates the points of
$V_{0,\gamma}$, the minimal polynomial $m_{\widetilde{u}}\in\Q[Y]$
of $\widetilde{u}$ in $V_{0,\gamma}$ and the parametrizations of
$Z_1\klk Z_n$ by the zeros of $m_{\widetilde{u}}$.

For this purpose, assume that we are given a linear form $\widetilde
u:=\widetilde u_1 Z_1\plp \widetilde u_n Z_n\in\Q[Z_1\klk Z_n]$
which separates the points of $V_{0,\gamma}$. Observe that the fact
that $\widetilde u$ is a separating linear form for $V_{0,\gamma}$
implies that $\widetilde u_i\not=0$ holds for $i=1\klk n$. Let
$Y,\widetilde{Y}$ be new indeterminates and let $m_1\klk
m_n\in\Q[Y]$ be the sequence of polynomials defined recursively by:
\begin{equation}\label{eq: resultant for diagonal syst}
m_1:={\widetilde u_1}^{-r_1}Y^{r_1}-q_1,\ m_i:=Res_{\widetilde{Y}}
\big(\widetilde u_i^{-r_i}(Y-\widetilde{Y})^{r_i}-q_i,
m_{i-1}(\widetilde{Y})\big)\ \mathrm{for}\ 2\le i\le n.
\end{equation}
We claim that the polynomial $m_n$ equals (up to scaling by a
nonzero element of $\Q$) the minimal polynomial $m_{\widetilde
u}\in\Q[Y]$ of the coordinate function induced by $\widetilde u$
in the $\Q$--algebra extension
$\Q\hookrightarrow\Q[V_{0,\gamma}]$. Indeed, for every $2\le i \le
n$, the polynomial $m_i(Y)$ is a linear combination of $\widetilde
u_i^{-r_i}(Y-\widetilde{Y})^{r_i}-q_i$ and
$m_{i-1}(\widetilde{Y})$ over $\Q[Y, \widetilde Y]$. Let
$u^{(i)}:=\widetilde u_1 Z_1\plp \widetilde u_iZ_i$ for $1\le i\le
n$. Then, the identity $\widetilde
u_i^{-r_i}(u^{(i)}-u^{(i-1)})^{r_i}-q_i=0$ holds in
$\Q[V_{0,\gamma}]$. Thus, assuming inductively that
$m_{i-1}(u^{(i-1)})=0$ in $\Q[V_{0,\gamma}]$, it follows that
$m_i(u^{(i)})=0$ in $\Q[V_{0,\gamma}]$ as well. Taking into
account that $\deg m_n\le r_1\dots r_n$ and that $m_{\widetilde
u}$ is a nonzero polynomial of degree $D_\gamma:=r_1\cdots
r_n=\#(V_{0,\gamma})$, we conclude that our claim holds.

In order to compute the polynomial $m_{\widetilde{u}}$, we compute
the resultants in (\ref{eq: resultant for diagonal syst}). Since
the resultant $Res_{\widetilde{Y}} \big(\widetilde
u_i^{-r_i}(Y-\widetilde{Y})^{r_i}-q_i,
m_{i-1}(\widetilde{Y})\big)$ is a polynomial of $\Q[Y]$ of degree
$r_1\cdots r_i$, by univariate interpolation in the variable
$\widetilde{Y}$ we reduce its computation to the computation of
$r_1\cdots r_i+1$ resultants of univariate polynomials in
$\Q[\widetilde Y]$. This interpolation step requires $O\big({\sf
M}(r_1^2\cdots r_i^2)\big)$ arithmetic operations in $\Q$ and does
not require any division by a nonconstant polynomial in the
coefficients $\widetilde{u}_1\klk \widetilde{u}_n$ (see, e.g.,
\cite{BoLeSc03}, \cite{BoSc05}). Each univariate resultant can be
computed using the algorithms in e.g. \cite{BiPa94}, \cite{GaGe99}
with ${\sf M}(r_1\cdots r_i)$ arithmetic operations in $\Q$.
Altogether, we obtain an algorithm for computing the minimal
polynomial $m_{\widetilde u}$ which performs $O\big({\sf
M}(D_\gamma^2)\big)$ arithmetic operations in $\Q$.

Next, we extend this algorithm to an algorithm for computing a
geometric solution of $V_{0,\gamma}$ as explained in Subsection
\ref{subsect: geometric sol}. We obtain the following result:
\begin{proposition}
\label{prop: first step computation geo sol V_0} Suppose that the
coefficients of the linear form $\widetilde u$ are randomly chosen
in the set $\{1\klk 4n\rho D_\gamma^3\}$, where $\rho$ is a fixed
positive integer. Then the algorithm described above computes a
geometric solution of the variety $V_{0,\gamma}$ (in the
coordinate system $Z_1,\dots, Z_n$) with error probability at most
$1/\rho$ using $O\big(n{\sf M}( D_\gamma^2)\big)$ arithmetic
operations in $\Q$.
\end{proposition}
\begin{proof}
As proved by our previous arguments, it is clear that the
algorithm described computes a geometric solution of
$V_{0,\gamma}$ with the stated number of arithmetic operations in
$\Q$. There remains to analyze its error probability.

The only probabilistic step of the algorithm is the choice of the
coefficients of the linear form $\widetilde u$, which must satisfy
two requirements. First, $\widetilde u$ must separate the points
of the variety $V_{0, \gamma}$. Since $V_{0,\gamma}$ consists of
$D_\gamma$ distinct points of $\A^n$, from Theorem \ref{th:
Zippel-Schwartz} it follows that for a random choice of the
coefficients of $\widetilde u$ in the set $\{1,\dots,4n\rho
D_\gamma^3\}$, the linear form $\widetilde u$ separates the points
of $V_{0,\gamma}$ with error probability at most $1/4n\rho
D_\gamma\le 1/2\rho$.

The second requirement concerns the computation of the univariate
resultants of the generic versions of the polynomials in (\ref{eq:
resultant for diagonal syst}). This is required in order to extend
the algorithm for computing the minimal polynomial
$m_{\widetilde{u}}$ to an algorithm for computing a geometric
solution of the variety $V_{0,\gamma}$. We use a fast algorithm
for computing resultants over $\Q(\Lambda)$ based on the Extended
Euclidean Algorithm (EEA for short). We shall perform the EEA over
the ring of power series $\Q[\![\Lambda-\widetilde u]\!]$,
truncating all the intermediate results up to order 2. Therefore,
the choice of the coefficients of $\widetilde u$ must guarantee
that all the elements of $\Q[\Lambda]$ which have to be inverted
during the execution of the EEA are invertible elements of the
ring $\Q[\![\Lambda-\widetilde u]\!]$.

For this purpose, we observe that, similarly to the proof of
\cite[Theorem 6.52]{GaGe99}, one deduces that all the denominators
of the elements of $\Q(\Lambda)$ arising during the application of
the EEA to the generic version of the polynomials $\widetilde
u_i^{-r_i}(\alpha-u^{(i-1)})^{r_i}-q_i$ and $m_{i-1}(u^{(i-1)})$
are divisors of at most $r_1\cdots r_{i-1}$ polynomials of
$\Q[\Lambda]$ of degree $2r_1\cdots r_i$ for any $\alpha\in\Q$.
This EEA step must be executed for $r_1\cdots r_i$ distinct values
of $\alpha\in\Q$, in order to perform the interpolation step.
Hence the product of the denominators arising during all the
applications of the EEA has degree at most $2nD_\gamma^3$.
Therefore, from Theorem \ref{th: Zippel-Schwartz} we conclude that
for a random choice of its coefficients in the set
$\{1,\dots,4n\rho D_\gamma^3\}$, the linear form $\widetilde u$
satisfies our second requirement with error probability at most
$1/2\rho$.

The lemma follows putting both error probability estimates
together.
\end{proof}

Finally, we compute a geometric solution of the variety
$V_{0,\gamma}$ in the original coordinate system defined by
$X_1,\dots, X_n$.

For this purpose, we compute the minimal polynomial $m_u\in\Q[Y]$
of a linear form $u = u_1 X_1 + \cdots + u_n X_n \in\Q[\xon]$ in
$V_{0,\gamma}$. Let $V_{0,\gamma}:=\{\mathrm{x}_0^{(1,\gamma)}\klk
\mathrm{x}_0^{(D_\gamma,\gamma)}\}$. Then we have
$m_u(Y)=\prod_{j=1}^{D_\gamma}(Y-u(\mathrm{x}_0^{(j,\gamma)}))$.
In order to compute $m_u$, we use the polynomials
$m_{\widetilde{u}},\widetilde{w}_1\klk \widetilde{w}_n$ which form
the previously computed geometric solution of $V_{0,\gamma}$ in
the variables $Z_1,\dots, Z_n$: from the identities
$X_i:=Z_1^{k_{1,i}^{(\gamma)}}\cdots Z_n^{k_{n,i}^{(\gamma)}}$
$(1\le i\le n)$ we deduce that $m_u$ equals the minimal polynomial
of the image of the projection $\eta_u: V_{0,\gamma}\to\A^1$
defined by $\eta_u^{(\gamma)}(z_1\klk z_n):=
\sum_{i=1}^nu_iz_1\!\!{}^{k_{1,i}^{(\gamma)}}\cdots
z_n\!\!{}^{k_{n,i}^{(\gamma)}}$. Now, the identities
$Z_i=\widetilde w_i(\widetilde u)$, which hold in
$\Q[V_{0,\gamma}]$ for $1\le i\le n$, imply that
\begin{equation}\label{eq: computation minimal V_0 bis}
u=\sum_{i=1}^nu_i\big(\widetilde w_1(\widetilde u)
\big)^{k_{1,i}^{(\gamma)}} \cdots \big (\widetilde w_n(\widetilde
u)\big)^{k_{n,i}^{(\gamma)}}\end{equation}
holds in $\Q[V_{0,\gamma}]$, from which we easily conclude that
$m_u$ satisfies the following identity:
\begin{equation}\label{eq: computation minimal V_0}
m_u(Y)=Res_{\widetilde{Y}}\Big(Y-\sum_{i=1}^nu_i \big(\widetilde
w_1(\widetilde{Y})\big)^{k_{1,i}^{(\gamma)}} \cdots \big
(\widetilde
w_n(\widetilde{Y})\big)^{k_{n,i}^{(\gamma)}},m_{\widetilde
u}(\widetilde{Y})\Big).
\end{equation}
We compute the monomials $\big(\widetilde w_1(\widetilde
u)\big){}^{k_{1,i}^{(\gamma)}} \cdots \big (\widetilde
w_n(\widetilde u)\big){}^{k_{n,i}^{(\gamma)}}$ $(1\le i\le n)$ in
the right--hand side of (\ref{eq: computation minimal V_0 bis})
modulo $m_{\widetilde u}(Y)$, with
$O\big(n^2\log(\max_{i,j}|k_{i,j}^{(\gamma)}|){\sf
M}(D_\gamma)\big)$ additional arithmetic operations in $\Q$. From
(\ref{eq: estimate height Smith form}) it follows that
$$O\big(n^2\log(\max_{i,j}|k_{i,j}^{(\gamma)}|){\sf
M}(D_\gamma)\big)=O\big(n^3\log(n\|\mathcal{E}_\gamma\|){\sf
M}(D_\gamma)\big),$$
where $\mathcal{E}_\gamma$ is the matrix of the exponents of the
cell corresponding to the inner normal $\gamma$. Observe that all
these steps are independent of the coefficients of the linear form
$u$ we are considering and therefore do not introduce any division
by a nonconstant polynomial in the coefficients $u_1\klk u_n$.

In the next step we compute the right--hand side of (\ref{eq:
computation minimal V_0 bis}) modulo $m_{\widetilde{u}}(Y)$, with
$O\big(nD_\gamma\big)$ arithmetic operations in $\Q$. Then we
compute the resultant (\ref{eq: computation minimal V_0}) by a
process which interpolates (\ref{eq: computation minimal V_0}) in
the variable $Y$ to reduce the question to the computation of
$D_\gamma+1$ univariate resultants, in the same way as for the
computation of the resultants in (\ref{eq: resultant for diagonal
syst}). This requires $O\big({\sf M}(D_\gamma)^2\big)$ arithmetic
operations in $\Q$.

If the linear form $u$ separates the points of $V_{0,\gamma}$,
then we can extend the algorithm for computing $m_u(Y)$ to an
algorithm for computing a geometric solution of $V_{0,\gamma}$
with the algorithm underlying the proof of Lemma \ref{lemma:
min/resgeom}. This extension requires that the coefficients
$u_1\klk u_n$ of the linear form $u$ do not annihilate the
denominators in $\Q[\Lambda]$ which arise from the application of
the algorithm described above to the generic version $\Lambda_1
X_1\plp\Lambda_nX_n$ of the linear form $u$. Such denominators
arise only during the computation of the generic version of the
resultant (\ref{eq: computation minimal V_0}). Hence, with a
similar analysis as in the proof of Proposition \ref{prop: first
step computation geo sol V_0}, we conclude that, if the
coefficients of $u$ are chosen randomly in the set $\{1\klk 4\rho
D_\gamma^3\}$, then the error probability of our algorithm is
bounded by $1/\rho$. In conclusion, we have:
\begin{proposition}\label{prop: second step computation geo sol V_0 gamma}
Suppose that we are given a geometric solution of $V_{0,\gamma}$
in the coordinate system $Z_1\klk Z_n$, as provided by the
algorithm underlying Proposition \ref{prop: first step computation
geo sol V_0}, and the coefficients of the linear form $u$ are
randomly chosen in the set $\{1\klk 4\rho D_\gamma^3\}$, where
$\rho$ is a fixed positive integer. Then the algorithm described
above computes a geometric solution of the variety $V_{0,\gamma}$
with error probability at most $1/\rho$ using
$O\big(n^3\log(n\|\mathcal{E}_\gamma\|){\sf M}(D_\gamma)^2\big)$
arithmetic operations in $\Q$.
\end{proposition}

Finally, from Propositions \ref{prop: first step computation geo
sol V_0} and \ref{prop: second step computation geo sol V_0 gamma}
and the fact that $\|E_\gamma\|$ is bounded by
$\mathcal{Q}:=2\max_{1\le i\le n}\{\|q\|;q\in\Delta_i\}$, we
immediately deduce the following result:
\begin{theorem}\label{th: computation geo sol V_0}
Suppose that the coefficients of the linear forms $\widetilde{u}$
and $u$ of the statement of Propositions \ref{prop: first step
computation geo sol V_0} and \ref{prop: second step computation
geo sol V_0 gamma} are chosen at random in the set $\{1\klk 4n\rho
D^3\}$, where $\rho$ is a fixed positive integer. Then the
algorithm underlying Propositions \ref{prop: first step
computation geo sol V_0} and \ref{prop: second step computation
geo sol V_0 gamma} computes a geometric solution of the varieties
$V_{0,\gamma}$ for all $\gamma\in\Gamma$ with error probability at
most $2/\rho$ using $O\big(n^3\log(n\mathcal{Q}){\sf M}(D)^2\big)$
arithmetic operations in $\Q$.
\end{theorem}
%
%
\subsection{The computation of a geometric solution of
the first deformation} \label{subsect: comput geo sol V hat}
The second step of our algorithm consists in the computation of a
geometric solution of the curve $\widehat{V}$ of (\ref{eq: def V
hat}). This will be done by ``lifting'' the geometric solutions of
the varieties $V_{0,\gamma}$ computed in the previous section for
all $\gamma\in\Gamma$.

We recall the definition of the variety $\widehat{V}$. Let $I$
denote the ideal of $\Q[X,T]$ generated by the polynomials
$\widehat h_1, \dots, \widehat h_n$ of (\ref{eq: def h hat}),
which form the polyhedral deformation of the generic polynomials
$h_1\klk h_n$, and let $J$ denote the Jacobian determinant of
$\widehat{h}_1\klk \widehat{h}_n$ with respect to the variables
$\xon$. Let $V(I)$ be the set of common zeros in
$\mathbb{A}^{n+1}$ of $\widehat h_1, \dots, \widehat h_n$. Then
$\widehat{V}:=V(I:J^\infty)$.

Alternatively, let ${\pi}:V(I) \to \A^1$ be the linear projection
defined by ${\pi}(x,t)=t$. Consider the decomposition of $V(I)$
into its irreducible components $V(I)= \bigcup_{i=1}^{r+s}
\mathcal{C}_i$. Suppose that the restriction
${\pi}|_{\mathcal{C}_i}:\mathcal{C}_i\to\A^1$ of the projection
$\pi$ is dominant for $1\le i\le r$ and is not dominant for
$r+1\le i\le s$. We shall show that
$\widehat{V}:=\bigcup_{i=1}^r\mathcal{C}_i$ holds, i.e.,
$\widehat{V}$ is the union of all the irreducible components of
$V(I)$ which project dominantly over $\mathbb{A}^1$. Furthermore,
we shall show that $\widehat{V}\subset \A^{n+1}$ is a curve which
constitutes a suitable deformation of the variety defined by the
system $h_1=\dots=h_n=0$. For this purpose, we shall use the
following technical lemma:
\begin{lemma}\label{lemma: variedad dom}
Let $F_1,\dots, F_n \in \Q[X,T]$  and $\mathcal{V}:=\{ (x,t) \in
\A^{n+1} : F_1(x,t) = 0, \dots, F_n(x,t) = 0\}$. Set $I:=
(F_1,\dots, F_n) \subset \Q[X,T]$ and let $J$ denote the Jacobian
determinant of $F_1,\dots, F_n$ with respect to the variables $X$.
Consider the linear projection $\pi: \mathcal{V} \to \A^1$ defined
by $\pi(x,t):= t$.  Assume that $\#\pi^{-1}(t)\le D$ holds for
generic values of $t\in \A^1$ and that there exists a point $t_0
\in \A^1$ such that the fiber $\pi^{-1}(t_0)$ is a
zero-dimensional variety of degree $D$ with $J(x, t_0) \ne 0$ for
every $(x, t_0) \in \pi^{-1}(t_0)$.

Let $\mathcal{V}_{\rm dom}$ be the union of all the irreducible
components $\mathcal{C}$ of $\mathcal{V}$ with
$\overline{\pi(\mathcal{C})} = \A^1$. Then:
\begin{itemize}
\item $\mathcal{V}_{\rm dom}$ is a nonempty equidimensional
variety of dimension $1$. \item $\mathcal{V}_{\rm dom}$ is the
union of all the irreducible components of $\mathcal{V}$ having a
non-empty intersection with $\pi^{-1}(t_0)$. \item
$\mathcal{V}_{\rm dom} = V(I:J^\infty)$. \item The restriction
$\pi|_{\mathcal{V}_{\rm dom}} : \mathcal{V}_{\rm dom} \to \A^1$ is
a dominant map of degree $D$.
\end{itemize}
\end{lemma}
\begin{proof} First we observe that $\dim(\mathcal{C})\ge 1$
for each irreducible component $\mathcal{C}$ of $\mathcal{V}$,
since $\mathcal{V}$ is defined by $n$ polynomials in an
$(n+1)$-dimensional space.

Let $\mathcal{C}$ be an irreducible component of $\mathcal{V}$ for
which $\pi^{-1}(t_0) \cap \mathcal{C} \ne \emptyset$ holds.
Consider the restriction $\pi|_\mathcal{C} : \mathcal{C} \to \A^1$
of the projection map $\pi$. Then we have that
$\pi|_\mathcal{C}^{-1}(t_0)$ is a nonempty zero-dimensional
variety, which implies that the generic fiber of
$\pi|_\mathcal{C}$ is either zero-dimensional or empty. Since
$\dim(\mathcal{C}) \ge 1$, the Theorem on the Dimension of Fibers
implies that $\dim(\mathcal{C}) = 1$ and that $\pi|_\mathcal{C}:
\mathcal{C} \to \A^1$ is a dominant map with generically-finite
fibers. This shows that $\mathcal{C} \subset \mathcal{V}_{\rm
dom}$ and, in particular, that $\mathcal{V}_{\rm dom}$ is
nonempty.

Conversely, we have that $\pi^{-1}(t_0) \cap \mathcal{C}\ne
\emptyset$ holds for any irreducible component $\mathcal{C}$ of
$\mathcal{V}_{\rm dom}$. Indeed, assume on the contrary the
existence of an irreducible component $\mathcal{C}_0$ not
satisfying this condition. Then, there is a point $t_1 \in \A^1$
having a finite fiber $\pi^{-1}(t_1)$ such that
$\pi|_{\mathcal{C}_0}^{-1}(t_1)$ and
$\pi|_{\mathcal{C}}^{-1}(t_1)$ have maximal cardinality for every
$\mathcal{C}$ with $\mathcal{C} \cap \pi^{-1}(t_0) \ne \emptyset$.
This implies that $\# \pi^{-1}(t_1) > \# \pi^{-1}(t_0)=D$, leading
to a contradiction.

We conclude that $\mathcal{V}_{\rm dom}$ is the nonempty
equidimensional variety of dimension $1$ which consists of all the
irreducible components $\mathcal{C}$ of $\mathcal{V}$ with
$\pi^{-1}(t_0) \cap \mathcal{C} \ne \emptyset$. Furthermore, this
shows that the restriction $\pi|_{\mathcal{V}_{\rm dom}} :
\mathcal{V}_{\rm dom} \to \A^1$ is a dominant map of degree $D$.

Finally we show that the identity $\mathcal{V}_{\rm dom} = V(I:
J^\infty)$ holds.   First, note that the irreducible components of
$V(I:J^\infty)$ are all the irreducible components of
$\mathcal{V}$ where the Jacobian $J$ does not vanish identically.
Thus, it is clear that $\mathcal{V}_{\rm dom} \subset
V(I:J^\infty)$, since $J$ does not vanish at the points of
$\pi^{-1}(t_0) \cap \mathcal{C}$ for each irreducible component
$\mathcal{C}$ of $\mathcal{V}_{\rm dom}$.   On the other hand, if
$\mathcal{C}$ is an irreducible component of $\mathcal{V}$ for
which the projection $\pi|_\mathcal{C}: \mathcal{C} \to \A^1$ is
not dominant, then $\mathcal{C}$ is the set of common zeros of the
polynomials $F_1,\dots, F_n,T - t_\mathcal{C}$ for some value
$t_\mathcal{C}$. Since $\dim(\mathcal{C}) \ge 1$, we have that the
Jacobian matrix $\partial(F_1,\dots, F_n, T-
t_\mathcal{C})/\partial(X_1,\dots, X_n, T)$ is singular at every
point $(x, t_\mathcal{C})$ of $\mathcal{C}$. Hence, its
determinant, which equals $J$, vanishes over $\mathcal{C}$.
\end{proof}

Now we return to the study of the variety $\widehat V$ and show
that the assumptions of Lemma \ref{lemma: variedad dom} hold.
Observe that $\pi^{-1}(t) = V_t \times \{ t\}$ holds for every
$t\in \A^1$, where $V_t:=\{ x\in \A^n: \widehat h_1(x,t) =
0,\dots, \widehat h_n(x,t) = 0\}$. Furthermore, the polynomials
$\widehat h_1(X,t), \dots, \widehat h_n(X,t)$ are obtained by a
suitable substitution of the variables $\Omega$ of the generic
polynomials $H_1,\dots, H_n\in \Q[\Omega, X]$ with supports
$\Delta_1,\dots, \Delta_n$ introduced in (\ref{eq: His}). Indeed,
if $c = (c_1,\dots, c_n)$ is the vector of coefficients of
$h_1,\dots, h_n$, the coefficient vector of $\widehat h_i(X,t)$
$(1\le i\le n)$ is $(c_{i,q} t^{\omega_i(q)})_{q\in \Delta_i}$ for
every $t\in \A^1$. By Lemma \ref{lemma: sparse system is 0dim},
there exists a nonzero polynomial $P^{(0)} \in \Q[\Omega]$ such
that for any $c'= (c_1',\dots,c_n')$ with $P^{(0)}(c') \not=0$ the
associated sparse system defines a zero-dimensional variety. In
particular, the coefficients $c = (c_1,\dots, c_n)$ of our input
polynomials $h_1:= H_1(c_1,X),\dots, h_n = H_n(c_n,X)$ satisfy
$P^{(0)}(c) \ne 0$. This shows that the polynomial $P^{(0)}_T\in
\Q[T]$ obtained by substituting $\Omega_{i,q} \mapsto c_{i,q}
T^{\omega_i(q)}$ $(1\le i \le n, \, q\in \Delta_i)$ in the
polynomial $P^{(0)}$ is nonzero, since it does not vanish at $T =
1$. We conclude that $V_t$ is a zero-dimensional variety for all
but a finite number of $t\in \A^1$. Thus, $\pi^{-1}(t)$ is finite
for generic values of $t\in \A^1$.

Finally, by condition $\mathsf{(H1)}$, the fiber
$\pi^{-1}(1)=V(h_1,\dots, h_n) \times \{ 1\}$ is a
zero--dimensional variety of degree $D = \deg(\pi)$ and the
Jacobian determinant $J:=\det(\partial\widehat h_i/\partial
X_j)_{1\le i,j\le n}$ does not vanish at any of its points. On the
other hand, the fact that $\#\pi^{-1}(t)\le D$ holds for generic
values $t\in\A^1$ follows from the BKK theorem.

This shows that the variety $V(I)$ and its defining polynomials
$\widehat h_1,\dots, \widehat h_n$ satisfy all the assumptions of
Lemma \ref{lemma: variedad dom}. Thus, we have:

\begin{lemma}\label{lemma: V hat es curva}
The variety $\widehat{V}\subset \A^{n+1}$ is a curve. Furthermore,
every irreducible component of $\widehat{V}$ has a nonempty
intersection with the fiber $\pi^{-1}(1)$ of the projection map
$\pi:\widehat V\to\A^1$.
\end{lemma}
%
%
\subsubsection{Generic linear projections of $\widehat{V}$.}
In order to compute a geometric solution of the space curve
$\widehat{V}$, we shall first exhibit a procedure for computing
the minimal polynomial of a generic linear projection of
$\widehat{V}$. Let $u\in\Q[\xon]$ be a linear form which separates
the points of the ``initial varieties'' $V_{0,\gamma}$ for all the
inner normals $\gamma:=(\gamma_1\klk \gamma_{n+1})$ of the lower
facets of the polyhedral deformation under consideration. Let
$\pi_u:\widehat{V}\to\A^2$ be the morphism defined by
$\pi_u(x,t):=(t,u(x))$. Since the projection map
$\pi:\widehat{V}\to\A^1$ defined by $\pi(x,t):=t$ is dominant, it
follows that the Zariski closure of the image of $\pi_u$ is a
$\Q$--definable hypersurface of $\A^{2}$. Denote by
$M_u\in\Q[T,Y]$ a minimal defining polynomial for this
hypersurface. For the sake of the argument, we shall assume
further that the identity $\deg(\pi)=D$, and thus $\deg_YM_u=D$,
hold.

We can apply estimate (\ref{eq: estimate height polyhedral def})
of Lemma \ref{lemma: estimate sparse height} in order to estimate
$\deg_T M_u$ in combinatorial terms (compare with \cite[Theorem
1.1]{PhSo06}). Indeed, let $\widehat Q_1,\dots, \widehat
Q_n\subset \R^{n+1}$ be the Newton polytopes of the polynomials
$\widehat{h}_1,\dots,\widehat{h}_n$ of (\ref{eq: def h hat}), and
let $\Delta\subset \R^{n+1}$ be the standard unitary simplex in
the plane $\{T = 0\}$. Then the following estimate holds:
\begin{equation}\label{eq: estimate height polyhedral def bis}
\deg_T M_u\le E:=MV_{n+1}(\Delta, \widehat Q_1,\dots, \widehat
Q_n).
\end{equation}
Furthermore, equality holds in (\ref{eq: estimate height
polyhedral def bis}) for a generic choice of the coefficients of
the polynomials $\widehat{h}_{i}$ and the linear form $u$.

More precisely, we shall exhibit a procedure for computing the
unique monic multiple in $\Q(T)[Y]$ of $M_u$ of degree $D$. This
polynomial can be alternatively defined as explained in what
follows:

Since the projection map $\pi: \widehat{V}\to \A^1$ is dominant,
it induces an extension $\Q[T]\hookrightarrow \Q[\widehat{V}]$,
where $\Q[\widehat{V}]$ denotes the coordinate ring of
$\widehat{V}$. This variety being a curve, $\Q[\widehat{V}]$ turns
out to be a finitely generated $\Q[T]$-module. Thus, tensoring
with $\Q(T)$, we deduce that $\Q[\widehat{V}]\otimes \Q(T)$ is a
$\Q(T)$--vector space of finite dimension. We claim that
$\Q[\widehat{V}]\otimes \Q(T) = \Q[V(I)]\otimes \Q(T)$ holds.
Indeed, since $\widehat V$ consists of the irreducible components
of $V(I)$ which are mapped dominantly onto $\A^1$ by the
projection $\pi$, for each of the remaining irreducible components
$\mathcal{C}$ of $V(I)$, the set $\pi(\mathcal{C})\subset\C$ is a
zero--dimensional $\Q$--definable variety. This implies that
$I(\mathcal{C}) \cap \Q[T]\not=\{0\}$ holds.

Let $\widehat m_u$ be the minimal polynomial of $u$ in the
extension $\Q(T) \hookrightarrow \Q[\widehat{V}]\otimes \Q(T)$.
The fact that $\Q[\widehat{V}] \otimes \Q(T)$ is
finite--dimensional $\Q(T)$--vector space shows that the affine
variety
$\mathbb{V}:=\{\bar{x}\in\A^n(\overline{\Q}(T)^*):\widehat{h}_1(\bar{x})=0\klk
\widehat{h}_n(\bar{x})=0\}$ has dimension zero. Here
$\overline{\Q}(T)^*:= \bigcup_{q \in
\N}\overline{\Q}(\!(T^{1/q})\!)$ denotes the field of Puiseux
series in the variable $T$ over $\overline{\Q}$ (see, e.g.,
\cite{Walker50}) and $\widehat{h}_1\klk\widehat{h}_n$ are
considered as elements of $\Q(T)[X]$. Our hypotheses imply that
there exist $D$ distinct $n$--tuples $x^{(\ell)}:=(x^{(\ell)}_1,
\ldots,x^{(\ell)}_n)\in(\overline{\Q}(T)^*)^n$ of Puiseux series
such that the following equalities hold in $\overline{\Q}(T)^*$
for $1\le\ell\le D$:
\begin{equation}\label{eq: puiseux series que anulan h hat}
\widehat{h}_1(x^{(\ell)},T)=0\:,\:\ldots\:,\: \widehat{h}_n
(x^{(\ell)},T)=0\end{equation}
(see \cite{HuSt95}). Since $\Q[\widehat{V}] \otimes \Q(T)$ is the
coordinate ring of the $\Q(T)$--variety $\mathbb{V}$, from
(\ref{eq: puiseux series que anulan h hat}) we deduce that the
dimension of $\Q[\widehat{V}] \otimes \Q(T)$ over $\Q(T)$ equals
$D$ . Moreover, since as a consequence of our assumptions
$\deg_Y\widehat{m}_u=D$ holds, we conclude that
\begin{equation}\label{eq: minimal u and puiseux}
\widehat{m}_u=\prod_{\ell=1}^D\big(Y-u(x^{(\ell)})\big).\end{equation}
Since $M_u(T,u(X))\in I(\widehat{V})$, it follows that $M_u(T,
u(X)) = 0$ holds in $\Q[\widehat{V}] \otimes \Q(T)$, from which we
conclude that $M_u$ is a multiple of $\widehat m_u$ by a factor in
$\Q(T)[Y]$. Taking into account that both are polynomials of
degree $D$ in the variable $Y$ and that $\widehat m_u$ is monic in
this variable, we deduce that $\widehat m_u$ is the quotient of
$M_u$ by its leading coefficient.
%
%
\subsubsection{A procedure for computing $\widehat{m}_u$.}
Now we exhibit a procedure for computing the minimal polynomial
$\widehat{m}_u$, which is based on the expression (\ref{eq:
minimal u and puiseux}) of $\widehat{m}_u$ in terms of the Puiseux
expansions (\ref{eq: def series puiseux solucion de h gamma}).
Then we will apply Lemma \ref{lemma: min/resgeom} to this
procedure in order to obtain an algorithm for computing a
geometric solution of the curve $\widehat{V}$.

With notations as in Section \ref{subsect: polyhedral deform}, let
$\Gamma\subset \Z^{n+1}$ be the set of primitive integer vectors
of the form $\gamma:= (\gamma_1,\dots, \gamma_n, \gamma_{n+1}) \in
\Z^{n+1}$ with $\gamma_{n+1} >0$ for which there is a cell
$C=(C^{(1)},\ldots,C^{(s)})$ of type $(k_1,\dots, k_s)$ of the
subdivision of $\mathcal{A}$ induced by $\omega$ such that
$\widehat C$ has inner normal $\gamma$. As asserted in Section
\ref{subsect: polyhedral deform}, if $\gamma\in\Gamma$ is the
inner normal of the lifting $\widehat{C}$ of a cell $C$ of type
$(k_1\klk k_s)$, there exist $D_\gamma:=k_1!\cdots
k_s!\cdot\mathrm{vol}(C)$ vectors of Puiseux series
$x^{(j,\gamma)}:=(x_1^{(j,\gamma)}\klk
x_n^{(j,\gamma)})\in\A^n({\Q}(T)^*)$ $(1\le j\le D_\gamma)$ of the
form
$$x^{(j,\gamma)}_i:= \sum_{m\ge 0}x^{(j,\gamma)}_{i,m}
T^{\frac{\gamma_i+m}{\gamma_{n+1}}}$$
satisfying (\ref{eq: puiseux series que anulan h hat}).
Considering the projection of the branches of $\widehat{V}$
parametrized by the $D_\gamma$ vectors of Puiseux series
$x^{(j,\gamma)}$ for each $\gamma\in\Gamma$, we obtain the
following element $m_{\gamma}$ of
$\Q(\!(T^{{1}/{\gamma_{n+1}}})\!)[Y]$:
\begin{equation}\label{eq: minimal u, gamma}
m_{\gamma}:=\prod_{j=1}^{D_\gamma}\big(Y-u(x^{(j,\gamma)})\big).
\end{equation}
{From} (\ref{eq: mixed volume and volumes}) we conclude that
(\ref{eq: minimal u and puiseux}) may be expressed in the
following way:
\begin{equation}\label{eq: minimal u, puiseux and cells}
\widehat{m}_u=\prod_{\gamma\in\Gamma}m_{\gamma}.\end{equation}

Since $\widehat{m}_u$ belongs to $\Q(T)[Y]$ and its primitive
multiple $M_u\in\Q[T,Y]$ satisfies the degree estimate $\deg_T
M_u\le E$, in order to compute the dense representation of
$\widehat{m}_u$ we shall compute the Puiseux expansions of the
coefficients of the factors
$m_{\gamma}\in\Q(\!(T^{{1}/{\gamma_{n+1}}})\!)[Y]$ of
$\widehat{m}_u$ truncated up to order $2E$. Using Pad\'e
approximation it is possible to recover the dense representation
of $\widehat{m}_u$ from this data.

Fix $\gamma \in \Gamma$ and set $\x_m^{(j,\gamma)}:=
(x_{1,m}^{(j,\gamma)}\klk x_{n,m}^{(j,\gamma)})$ for every $m\ge
0$ and $1\le j \le D_\gamma$. Since
$$\widehat{h}_i\Big(\sum_{m\ge
0}x^{(j,\gamma)}_{1,m} T^{\frac{\gamma_1+m}{\gamma_{n+1}}}\klk
\sum_{m\ge 0}x^{(j,\gamma)}_{n,m}
T^{\frac{\gamma_n+m}{\gamma_{n+1}}},T\Big)=0$$
holds for $1\le j\le D_\gamma$ and $1\le i\le n$, we have
\begin{equation}\label{eq: relation h hat and h gamma}
\begin{array}{rcl}
0&=&T^{-m_i}\widehat{h}_i\big(\sum_{m\ge 0}x^{(j,\gamma)}_{1,m}
T^{\gamma_1+m}\klk \sum_{m\ge 0}x^{(j,\gamma)}_{n,m}
T^{\gamma_n+m},T^{\gamma_{n+1}}\big)\\[2mm]
&=&T^{-m_i}\widehat{h}_i\big(T^{\gamma_1}\sum_{m\ge
0}x^{(j,\gamma)}_{1,m}T^m\klk T^{\gamma_n}\sum_{m\ge
0}x^{(j,\gamma)}_{n,m} T^m,T^{\gamma_{n+1}}\big)\\[2mm] &=&
h_{i,\gamma}\big(\sum_{m\ge 0}\x_m^{(j,\gamma)}T^m,T\big),
\end{array}\end{equation}
according to (\ref{eq: def h gamma}). Therefore the polynomial
$m_{\gamma}(T^{\gamma_{n+1}},Y)\in\Q(\!(T)\!)[Y]$ can be expressed
in terms of the power series solutions
$\sigma^{(j,\gamma)}:=(\sigma_1^{(j,\gamma)}\klk
\sigma_n^{(j,\gamma)}):=\sum_{m\ge 0}\x_m^{(j,\gamma)}T^m$ $(1\le
j\le D_\gamma)$ of $h_{1,\gamma}\klk h_{n,\gamma}$. Indeed, from
(\ref{eq: minimal u, gamma}) it follows that
\begin{equation}\label{eq: minimal u gamma as projection}
\begin{array}{rcl}
m_{\gamma}(T^{\gamma_{n+1}},Y)&=&\prod_{j=1}^{D_\gamma}
\big(Y-\sum_{i=1}^nu_i\sum_{m\ge
0}x_{i,m}^{(j,\gamma)}T^{\gamma_i+m}\big)
\\[2mm] &=&\prod_{j=1}^{D_\gamma}
\big(Y-\sum_{m\ge 0}\sum_{i=1}^nu_ix_{i,m}^{(j,\gamma)}
T^{\gamma_i}T^m\big) \\[2mm] &=& \prod_{j=1}^{D_\gamma}
\big(Y-\sum_{m\ge 0}u_\gamma(\x_m^{(j,\gamma)})T^m\big)
\\[2mm] &=& \prod_{j=1}^{D_\gamma}
\big(Y-u_\gamma(\sum_{m\ge 0}\x_m^{(j,\gamma)}T^m)\big)\ =:\
m_{u_\gamma}(T,Y),
\end{array}\end{equation}
where $u_\gamma:=\sum_{i=1}^nu_iT^{\gamma_i}X_i$. We conclude that
the Laurent polynomial \linebreak
$m_{\gamma}(T^{\gamma_{n+1}},\!Y)\in\Q(\!(T)\!)[Y]$ may be
considered as the minimal polynomial $m_{u_\gamma}(T,\!Y)$ of the
projection induced by $u_\gamma$ on the subvariety $V_\gamma$ of
$\A^n(\overline{\Q}(T)^*)$ consisting of the set of power series
$\{\sigma^{(1,\gamma)}\klk \sigma^{(D_\gamma,\gamma)}\}$. This
remark will allow us to compute a suitable approximation to the
Laurent polynomial $m_{\gamma} (T^{\gamma_{n+1}},Y)$ in
$\Q(\!(T)\!)[Y]$.

In order to describe this approximation, we introduce the following
terminology: for $G,\widetilde{G}\in\overline{\Q}(\!(T)\!)$ and
$s\in\Z$, we say that $\widetilde{G}$ {\sf approximates} $G$ {\sf
with precision} $s$ in $\overline{\Q}(\!(T)\!)$ if the Laurent
series $G-\widetilde{G}$ has order at least $s+1$ in $T$. We shall
use the notation $G\equiv\widetilde{G}$ mod $(T^{s+1})$.
Furthermore, if $G,\widetilde{G}$ are two elements of a polynomial
ring $\overline{\Q}(\!(T)\!)[Y]$, we say that $\widetilde{G}$
approximates $G$ with precision $s$ if every coefficient
$\widetilde{a}\in\overline{\Q}(\!(T)\!)$ of $\widetilde{G}$
approximates the corresponding coefficient $a\in\overline{\Q}
(\!(T)\!)$ of $G$ with precision $s$ (in the sense of the previous
definition).

\begin{proposition}\label{prop: approx minimal u gamma}
Fix $\gamma:=(\gamma_1\klk\gamma_n)\in \Gamma$ and assume that a
geometric solution of the variety $V_{0,\gamma}$ is given, as
provided by Theorem \ref{th: computation geo sol V_0}. Assume
further that the coefficients of the linear form $u$ of the given
geometric solution of $V_{0,\gamma}$ are randomly chosen in the
set $\{1\klk 4\rho D_\gamma^3\}$ for a given $\rho\in\N$. Then
there is an algorithm which computes an approximation to the
polynomial $m_{u_\gamma}\in \Q(\!(T)\!)[Y]$ with precision
$2E\gamma_{n+1}$. The procedure requires
$O\big((nL_\gamma+n^{\Omega}){\sf M}(D_\gamma)\big({\sf
M}(M_\gamma){\sf M}(D_\gamma)+E\gamma_{n+1}\big)\big)$ arithmetic
operations in $\Q$ , where
$M_\gamma:=\max\{\gamma_1\klk\gamma_n\}$ and $L_\gamma$ is the
number of arithmetic operations required to evaluate the
polynomials $h_{i,\gamma}$ of {\rm (\ref{eq: def h gamma})}, and
has error probability at most $2/\rho$.
\end{proposition}
\begin{proof} Let notations and assumptions be as before.
In order to compute the required approximation of the polynomial
$m_{u_\gamma}$ we first compute the corresponding approximation of
the polynomials that form a geometric solution of the variety
$V_\gamma:=\{\sigma^{(j,\gamma)};1\le j\le D_\gamma\}$. Observe
that
\begin{eqnarray*}
\{\sigma^{(j,\gamma)}(0);1\le j\le
D_\gamma\} &=& \{\mathrm{x}_0^{(j,\gamma)};1\le j\le D_\gamma\} \\[1mm]
   &=& V(h_{1,\gamma}^{(0)}\klk
h_{n,\gamma}^{(0)})\cap(\C^*)^n \\[1mm]
   &=& V(h_{1,\gamma}(X,0)\klk
h_{n,\gamma}(X,0))\cap(\C^*)^n=V_{0,\gamma}
\end{eqnarray*}
holds. Since $\det(\partial h_{i,\gamma}(X,0)/\partial X_k)_{1\le
i,k\le n}(x_0^{(j,\gamma)})\not=0$ holds for $1\le j\le D_\gamma$,
we may apply of the global Newton iterator of \cite{GiLeSa01} (see
also \cite{Schost03}) in order to ``lift'' the given geometric
solution of $V_{0,\gamma}$ to the geometric solution of the
variety $V_\gamma$ associated to the linear form $u\in\Q[X]$ with
any prescribed precision.

Denote $m_{u,\gamma}^{(0)},w_{u,1,\gamma}^{(0)}\klk
w_{u,n,\gamma}^{(0)}\in\Q[Y]$ the polynomials which form the given
geometric solution of $V_{0,\gamma}$, as provided by the algorithm
underlying Theorem \ref{th: computation geo sol V_0}. Recall that
$m_{u,\gamma}^{(0)}\big(u(\mathrm{x}_0^{(j)})\big)=0$ and
$(\mathrm{x}_0^{(j,\gamma)})_i=w_{u,i,\gamma}^{(0)}
\big(u(\mathrm{x}_0^{(j)})\big)$ holds for $1\le i\le n$ and $1\le
j \le D_\gamma$.
The global Newton iterator is a recursive procedure whose $k$th
step computes approximations
$m_{u,\gamma}^{(k)},w_{u,1,\gamma}^{(k)}\klk w_{u,
n,\gamma}^{(k)}\in\Q[T,Y]$ of the polynomials
$m_{u,\gamma},w_{u,1,\gamma}\klk w_{u,n,\gamma}$ which form the
geometric solution of $V_\gamma$ associated with the linear form
$u$ with precision $2^k$ for any $k\ge 0$.

Assume without loss of generality that $\gamma_i\ge 0$ and
$0=\min\{\gamma_1\klk \gamma_n\}$ hold for $1\le i\le n$. Indeed,
if there exists $\gamma_i<0$, setting
$\gamma_{i_0}:=\min\{\gamma_1\klk \gamma_n\}$ we have
\begin{eqnarray*}
T^{-\gamma_{i_0}D_\gamma}
m_{u,\gamma}(T^{\gamma_{n+1}},T^{\gamma_{i_0}}Y)&=&
\prod_{j=1}^{D_\gamma}T^{-\gamma_{i_0}}
\Big(T^{\gamma_{i_0}}Y-\sum_{i=1}^nu_i\sum_{m\ge
0}x_{i,m}^{(j,\gamma)}T^{\gamma_i+m}\Big) \\
   &=& \prod_{j=1}^{D_\gamma}
\Big(Y-T^{-\gamma_{i_0}}\sum_{i=1}^nu_i\sum_{m\ge
0}x_{i,m}^{(j,\gamma)}T^{\gamma_i+m}\Big) \\
   &=& \prod_{j=1}^{D_\gamma}
\Big(Y-\sum_{i=1}^nu_i\sum_{m\ge
0}x_{i,m}^{(j,\gamma)}T^{\gamma_i-\gamma_{i_0}+m}\Big).
\end{eqnarray*}
 Since $\gamma_i-\gamma_{i_0}\ge 0$ holds for $1\le i\le n$,
this shows that the computation of an approximation
$m_{u_\gamma}:=m_{\gamma}(T^{\gamma_{n+1}},Y)$ can be easily
reduced to a situation in which $\gamma_i\ge 0$ holds for $1\le
i\le n$.

Note that the global Newton iterator cannot be directly applied in
order to compute the geometric solution of
$\{\sigma^{(j,\gamma)};1\le j\le D_\gamma\}$ associated with the
linear form $u_\gamma\in\Q[T][X]$, because the coefficients of
$u_\gamma$ are nonconstant polynomials of $\Q[T]$. Indeed, two
critical problems arise:
\begin{enumerate}
    \item\label{item: primer problema Newton Lecerf}
    Although by hypothesis $u_\gamma$ separates the points of
$V_\gamma$, it might not separate the points of $V_{0,\gamma}$ and
it is not clear from which precision on, the corresponding
approximations of the points of $V_\gamma$ are separated by
$u_\gamma$. Requiring $u_\gamma$ to be a separating form for all
the approximations of the points of $V_\gamma$ is an essential
hypothesis for the iterator of \cite{GiLeSa01} which cannot be
suppressed without causing a significant growth of the complexity
of the procedure (see \cite{Lecerf02}, \cite{Lecerf03}).
    \item\label{item: segundo problema Newton Lecerf}
    The iterator of \cite{GiLeSa01} makes critical use of the fact
that the coefficients of the linear form under consideration are
elements of $\Q$ in order to determine how a given precision can
be achieved.
\end{enumerate}
Nevertheless, we shall exhibit a modification of the procedure
which computes an approximation of $m_{u_\gamma}(T,Y)$ with
precision $2\gamma_{n+1}E$ without changing the asymptotic number
of arithmetic operations performed.

In order to circumvent (\ref{item: primer problema Newton Lecerf})
we require an additional generic condition to be satisfied by the
coefficients $u_1\klk u_n$ defining $u_\gamma:=\sum_{i=1}^n
u_iT^{\gamma_i}X_i$. Recall that $u_\gamma(\sigma^{(j,\gamma)})  =
\sum_{m\ge 0} \big(\sum_{i=1}^n u_i
x_{i,m-\gamma_i}^{(j,\gamma)}\big)T^{m}$ for every $1\le j \le
D_\gamma$, where $x_{i,m-\gamma_i}^{(j,\gamma)}:=0$ for $m<
\gamma_i$. To state this condition, we need the following claim:
\begin{claim} Set $M_\gamma:= \max
\{\gamma_1,\dots,\!\gamma_n\}$ and let $\Lambda_1,\dots,\Lambda_n$
be indeterminates over $\C[T\!,\!X]$. Then the following
inequality holds for every $1\le j,h\le D_\gamma$ with $j\not=h$:
$$\sum_{m=0}^{M_\gamma}\Big(\sum_{i=1}^n\Lambda_i\,
x_{i,m-\gamma_i}^{(j,\gamma)}\Big)T^{m}\not=
\sum_{m=0}^{M_\gamma}\Big(\sum_{i=1}^n\Lambda_i\,
x_{i,m-\gamma_i}^{(h,\gamma)}\Big)T^{m}.
$$
\end{claim}

\noindent {\it Proof of Claim.} Suppose on the contrary that there
exist $j\not= h$ such that \linebreak
$\sum_{m=0}^{M_\gamma}\big(\!\sum_{i=1}^n\Lambda_i\,
x_{i,m-\gamma_i}^{(j,\gamma)}\big)T^{m}=
\sum_{m=0}^{M_\gamma}\big(\!\sum_{i=1}^n\Lambda_i\,
x_{i,m-\gamma_i}^{(h,\gamma)}\big)T^{m}$. Substituting
$T^{-\gamma_i}\Lambda_i$ for $\Lambda_i$ in this identity for
$i=1,\dots,n$, we have $\sum_{m=0}^{M_\gamma}\sum_{i=1}^n
\Lambda_i\, x_{i,m-\gamma_i}^{(j,\gamma)}T^{m-\gamma_i}=
\sum_{m=0}^{M_\gamma}\sum_{i=1}^n \Lambda_i\,
x_{i,m-\gamma_i}^{(h,\gamma)}T^{m-\gamma_i}$, that is
$$\sum_{i=1}^n \sum_{m=0}^{M_\gamma-\gamma_i} \Lambda_i\,
x_{i,m}^{(j,\gamma)}T^{m}= \sum_{i=1}^n
\sum_{m=0}^{M_\gamma-\gamma_i} \Lambda_i\,
x_{i,m}^{(h,\gamma)}T^{m}.$$
Substituting 0 for $T$ in this identity, we deduce that
$$\sum_{i=1}^n\Lambda_i x_{i,0}^{(j,\gamma)}=\sum_{i=1}^n\Lambda_i
x_{i,0}^{(h,\gamma)},$$
which contradicts the fact that the vectors $\x_0^{(j,\gamma)} =
(x_{1,0}^{(j,\gamma)}, \dots, x_{n,0}^{(j,\gamma)})$ $(1\le j \le
D_\gamma)$ are all distinct. This finishes the proof of the claim.
\medskip

By the claim we see that the polynomial
$\sum_{m=0}^{M_\gamma}\big(\sum_{i=1}^n\Lambda_i\,
(x_{i,m-\gamma_i}^{(j,\gamma)}-x_{i,m-\gamma_i}^{(h,\gamma)})\big)T^{m}$
of $\Q[\Lambda][T]$ is nonzero, and therefore has a nonzero
coefficient $a_{j,h}\in\C[\Lambda]$ for every $1\le j<h\le
D_\gamma$. Consider the polynomial $A_\gamma(\Lambda):=\prod_{1\le
j<h\le D_\gamma}a_{j,h}\in\C[\Lambda]$. Since $a_{j,h}$ has degree
1 for every $1\le j<h\le D_\gamma$, it follows that $A$ has degree
$\binom{D_\gamma}{2}$. Furthermore, for every $(u_1\klk
u_n)\in\C^n$ with $A_\gamma(u_1\klk u_n)\not=0$, the corresponding
polynomial $u_\gamma:=\sum_{i=1}^nu_iT^{\gamma_i} X_i$ separates
the initial terms
$\sum_{m=0}^{M_\gamma}\mathrm{x}_{m}^{(j,\gamma)}T^m$ of the power
series $\sigma^{(j,\gamma)}$ $(1\le j\le D_\gamma)$.

{From} Theorem \ref{th: Zippel-Schwartz} we see that for a random
choice of the coefficients $u_1\klk u_n$ in the set $\{1\klk \rho
D_\gamma^2\}$ the linear form $u_\gamma$ separates the first
$M_\gamma$ terms of the points of $V_\gamma$ with probability at
least $1-1/\rho$. From now on we assume that $u_\gamma$ satisfies
this requirement.

The algorithm proceeds in three steps. First, it computes a
suitable approximation to the geometric solution of $V_\gamma$
associated to the linear form $u:=\sum_{i=1}^nu_iX_i$ by means of
$\kappa_0:=\lceil\log (M_\gamma+1)\rceil$ steps of the global
Newton iterator of \cite{GiLeSa01}. This approximation is used in
order to obtain the corresponding approximation
$m_{u_\gamma}^{(\kappa_0)},w_{u_\gamma,1}^{(\kappa_0)}\klk
w_{u_\gamma, n}^{(\kappa_0)}$ of the polynomials that form the
geometric solution of $V_\gamma$ associated with $u_\gamma$.
Finally, we apply an adaptation of the global Newton iterator
which takes as input the polynomials of the previous step
$m_{u_\gamma}^{(\kappa_0)},w_{u_\gamma,1}^{(\kappa_0)}\klk
w_{u_\gamma, n}^{(\kappa_0)}$ and outputs the required
approximation to the polynomials $m_{u_\gamma},w_{u_\gamma,1}\klk
w_{u_\gamma, n}$ that form the geometric solution of $V_\gamma$
associated with $u_\gamma$.

Now we consider the three steps above in detail. The first step
takes as input the given geometric solution
$m_{u,\gamma}^{(0)},w_{u,1,\gamma}^{(0)}\klk w_{u,n,\gamma}^{(0)}$
of $V_{0,\gamma}$, and performs $\kappa_0:=\lceil\log
(M_\gamma+1)\rceil$ times the global Newton iterator of
\cite{GiLeSa01} to obtain polynomials
$m_{u,\gamma}^{(\kappa_0)},w_{u,1,\gamma}^{(\kappa_0)}\klk w_{u,
n,\gamma}^{(\kappa_0)}\in\Q[T,Y]$ such that the following
conditions hold:\smallskip

    $(i)_{u,\kappa_0}$\quad
    $\,\deg_Ym_{u,\gamma}^{(\kappa_0)}=D_\gamma$ and
    $\deg_Tm_{u,\gamma}^{(\kappa_0)}\le M_\gamma$,

    $(ii)_{u,\kappa_0}$\quad
    $\deg_Yw_{u, i,\gamma}^{(\kappa_0)}<D_\gamma$ and
     $\deg_Tw_{u,i,\gamma}^{(\kappa_0)}\le M_\gamma$
     for $1\le i\le n$,

    $(iii)_{u,\kappa_0}$\quad
    $\!m_{u,\gamma}^{(\kappa_0)}\equiv\prod_{j=1}^{D_\gamma}
    \big(Y-\varphi_{\kappa_0}^{(j,\gamma)}\big)$ mod
    $(T^{M_\gamma+1})$,

    $(iv)_{u,\kappa_0}$\quad
    $\sigma_i^{(j,\gamma)}\equiv w_{u,i,\gamma}^{(\kappa_0)}
    \big(T,\varphi_{\kappa_0}^{(j,\gamma)}\big)$ mod
    $(T^{M_\gamma+1})$ for $1\le i\le n$.
\smallskip

\noindent Here $\varphi_{\kappa_0}^{(j,\gamma)}$ is the Taylor
expansion of order $2^{\kappa_0}$ of the power series
$u(\sigma^{(j,\gamma)})$, that is,
$\varphi_{\kappa_0}^{(j,\gamma)}:=\sum_{m=0}^{2^{\kappa_0}}
u(\mathrm{x}_m^{(j,\gamma)})T^m$ for $1\le j\le D_\gamma$.
\medskip

\noindent According to \cite[Proposition 7]{GiLeSa01}, it follows
that this step requires performing
$O\big((nL_\gamma+n^{\Omega}){\sf M}(D_\gamma){\sf M}
(M_\gamma)\big)$ arithmetic operations in $\Q$, where $L_\gamma$
denotes the number of arithmetic operations in $\Q$ required to
evaluate the polynomials $h_{i,\gamma}$ of (\ref{eq: def h
gamma}). Furthermore, in view of the application of Lemma
\ref{lemma: min/resgeom} it is important to remark that this step
does not involve any division by a nonconstant polynomial in the
coefficients $u_1\klk u_n$.

Next we discuss the second step. In this step we obtain
approximations \linebreak $m_{u_\gamma}^{(\kappa_0)},
w_{u_\gamma,1}^{(\kappa_0)}\klk w_{u_\gamma, n}^{(\kappa_0)}$ of
the polynomials that form the geometric solution of $V_\gamma$
associated with $u_\gamma$ with precision $2^{\kappa_0}\ge
M_\gamma$, namely
\begin{itemize}
   \item $\deg_Ym_{u_\gamma}^{(\kappa_0)}=D_\gamma$ and
    $\deg_Tm_{u_\gamma}^{(\kappa_0)}\le M_\gamma$,
    \item $\deg_Yw_{u_\gamma,i}^{(\kappa_0)}<D_\gamma$ and
     $\deg_Tw_{u_\gamma,i}^{(\kappa_0)}\le M_\gamma$
     for $1\le i\le n$,
    \item $m_{u_\gamma}^{(\kappa_0)}\equiv\prod_{j=1}^{D_\gamma}
    \big(Y-\phi_{\kappa_0}^{(j,\gamma)}\big)$ mod
    $(T^{M_\gamma+1})$,
    \item $\sigma_i^{(j,\gamma)}\equiv w_{u_\gamma,i}^{(\kappa_0)}
    \big(T,\phi_{\kappa_0}^{(j,\gamma)}\big)$ mod
    $(T^{M_\gamma+1})$ for $1\le i\le n$.
\end{itemize}
Here $\phi_{\kappa_0}^{(j,\gamma)}$ is the Taylor expansion of
$\phi^{(j,\gamma)}:=u_\gamma(\sigma^{(j,\gamma)})$ of order
$2^{\kappa_0}$ for $1\le j\le D_\gamma$.

{From} conditions $(i)_{u,\kappa_0}\!$--$\,(iv)_{u,\kappa_0}$ and
the elementary properties of the resultant it is easy to see that
$m_{u_\gamma}^{(\kappa_0)}$ satisfies the following identity:
\begin{equation}\label{eq: comput resultant m u tilde}
m_{u_\gamma}^{(\kappa_0)}(Y)=
Res_{\widetilde{Y}}\Big(Y-\sum_{i=1}^n u_i T^{\gamma_i}
w_{u,i,\gamma}^{(\kappa_0)}(\widetilde{Y}),\,
m_{u,\gamma}^{(\kappa_0)}(\widetilde{Y})\Big).
\end{equation}
The resultant of the right--hand side is computed mod
$(T^{M_\gamma+1})$ by interpolation in the variable $Y$ to reduce
the problem to the computation of $D_\gamma$ resultants, as
explained in the computation of the resultant in (\ref{eq:
computation minimal V_0}). These $D_\gamma$ resultants involve two
polynomials of $\Q[T,\widetilde{Y}]$ of degree in $\widetilde{Y}$
bounded by $D_\gamma$ and are computed mod $(T^{M_\gamma+1})$.
Hence we deduce that this step requires $O\big({\sf
M}(D_\gamma)D_\gamma{\sf M}(M_\gamma)\big)$ arithmetic operations
in $\Q$.

We apply Lemma \ref{lemma: min/resgeom} in order to extend this
procedure to an algorithm computing $m_{u_\gamma}^{(\kappa_0)},
w_{u_\gamma,1}^{(\kappa_0)}\klk w_{u_\gamma, n}^{(\kappa_0)}$. For
this purpose, we observe that a similar argument as in the proof
of Proposition \ref{prop: first step computation geo sol V_0}
proves that the denominators in $\Q[\Lambda]$ which arise during
the computation of the $D_\gamma$ resultants required to compute
the minimal polynomial of the generic version
$\sum_{i=1}^n\Lambda_iT^{\gamma_i}X_i$ of the linear form
$u_\gamma$ are divisors of a polynomial of $\Q[\Lambda]$ of degree
at most $4D_\gamma^3$. Applying Theorem \ref{th: Zippel-Schwartz}
we see that for a random choice of the coefficients $u_1\klk u_n$
in the set $\{1\klk 4\rho D_\gamma^3\}$ none of these denominators
are annihilated with probability at least $1-1/\rho$.

Finally, we consider the third step of the algorithm. For
$\kappa_1:=\lceil\log(2\gamma_{n+1}E+1)\rceil$, we apply
$\kappa_1- \kappa_0$ times an adaptation of the global Newton
iterator of \cite{GiLeSa01} to the polynomials
$m_{u_\gamma}^{(\kappa_0)}, w_{u_\gamma,1}^{(\kappa_0)}\klk
w_{u_\gamma, n}^{(\kappa_0)}$ computed in the previous step. In
the $k$th iteration step, we compute polynomials
$m_{u_\gamma}^{(k)},w_{u_\gamma,1}^{(k)}\klk w_{u_\gamma,
n}^{(k)}$ satisfying:
\begin{itemize}
    \item $\deg_Ym_{u_\gamma}^{(k)}=D$ and
    $\deg_Tm_{u_\gamma}^{(k)}\le 2^k$,
    \item $m_{u_\gamma}^{(k)}=\prod_{j=1}^{D_\gamma}
    (Y-\phi_k^{(j,\gamma)})$,
    \item $\deg_Y w_{u_\gamma,i}^{(k)} <D$ and
     $\deg_Tw_{u_\gamma,1}^{(k)}\le 2^k$
     for $1\le i\le n$,
    \item $\sigma_i^{(j,\gamma)}\equiv w_{u_\gamma,i}^{(k)}
    (T,\phi_k^{(j,\gamma)})$ mod
    $(T^{2^k+1})$ for $1\le i\le n$.
\end{itemize}
Here $\phi_k^{(j,\gamma)}$ is the Taylor expansion of
$\phi^{(j,\gamma)}:=u_\gamma(\sigma^{(j,\gamma)})$ of order $2^k$
for $1\le j\le D_\gamma$.
In particular, it follows that $m_{u_\gamma}^{(\kappa_1)}$ is the
required approximation to $m_{u_\gamma}$ with precision
$2\gamma_{n+1}E$.

Fix $\kappa_0<k\le\kappa_1$. We briefly describe how we can obtain
an approximation with precision $2^k$ of the polynomials that form
the geometric solution of $V_\gamma$ associated to the linear form
$u_\gamma$ from an approximation with precision $2^{k-1}$.
Similarly to \cite{GiLeSa01}, set
$\Delta_k(T,Y):=u_\gamma(\widetilde{w}_{u_\gamma}^{(k)})-
u_\gamma(w_{u_\gamma}^{(k-1)})=u_\gamma(\widetilde{w}_{u_\gamma}^{(k)})-
Y$, where $\widetilde{w}_{u_\gamma}^{(k)}$ is the result of
applying a ``classical Newton step'' to ${w}_{u_\gamma}^{(k-1)}$,
as described in \cite{GiLeSa01}. Furthermore, write
$\Delta_m(T,Y):=T^{-1-2^{k-1}}(m_{u_\gamma}^{(k)}-
m_{u_\gamma}^{(k-1)})$. Since
$m_{u_\gamma}^{(k)}(Y+\Delta_k)\equiv
0\mod(T^{2^k+1},m_{u_\gamma}^{(k-1)})$ holds (see \cite[\S
4.2]{DuLe06}), it follows that
\begin{eqnarray*}0\equiv
m_{u_\gamma}^{(k)}(Y\!\!+\Delta_k)\!\!\!\!&\equiv&\!\!\!\!
m_{u_\gamma}^{(k-1)}(Y\!\!+\Delta_k)+T^{2^{k-1}+1}\Delta_m(Y\!\!+\Delta_k)
\mod(T^{2^k+1},m_{u_\gamma}^{(k-1)})\\
\!\!\!\!&\equiv&\!\!\!\! \Delta_k\frac{\partial
m_{u_\gamma}^{(k-1)}\!\!\!\!}{\partial
Y}(Y)+T^{2^{k-1}+1}\Delta_m(Y)
\mod(T^{2^k+1},m_{u_\gamma}^{(k-1)}).
\end{eqnarray*}
We conclude that the following congruence relation holds:
\begin{equation}\label{eq: congruence relation adaptation Newton Lecerf}
\quad\qquad m_{u_\gamma}^{(k)}\equiv
m_{u_\gamma}^{(k-1)}-\Big(\Delta_k\frac{\partial
m_{u_\gamma}^{(k-1)}\!\!\!\!}{\partial Y}\mod
m_{u_\gamma}^{(k-1)}\Big)\mod(T^{2^k+1}).
\end{equation}

A similar argument proves the following congruence relation
\begin{equation}
\label{eq: congruence relation adaptation Newton Lecerf 2}
w_{u_\gamma,i}^{(k)}\equiv
\widetilde{w}_{u_\gamma,i}^{(k-1)}-\Big(\Delta_k\frac{\partial
\widetilde{w}_{u_\gamma,i}^{(k-1)}\!\!\!\!}{\partial Y}\!\!\mod
m_{u_\gamma}^{(k-1)}\Big)\!\!\mod(T^{2^k+1})\ \mathrm{for}\ 1\le
i\le n.
\end{equation}

Each iteration of our adaptation of the global Newton iteration is
based on (\ref{eq: congruence relation adaptation Newton Lecerf})
and (\ref{eq: congruence relation adaptation Newton Lecerf 2}),
which are extensions of the corresponding congruence relations of
\cite{GiLeSa01}. We first compute $\widetilde{w}_{u_\gamma}^{(k)}$
by a standard Newton--Hensel lifting, and then evaluate the
expressions (\ref{eq: congruence relation adaptation Newton
Lecerf}) and (\ref{eq: congruence relation adaptation Newton
Lecerf 2}). With a similar analysis as in \cite[Proposition
7]{GiLeSa01} we conclude that the whole procedure requires
$O\big((nL_\gamma+n^{\Omega}){\sf M}(D_\gamma)E\gamma_{n+1}\big)$
arithmetic operations in $\Q$.

Finally, combining the complexity estimates of the three steps
above and the probability of achievement of the two generic
conditions imposed to the coefficients $u_1\klk u_n$, we deduce
the statement of the proposition.
\end{proof}

Using the algorithm of the statement of Proposition \ref{prop:
approx minimal u gamma} for all $\gamma\in\Gamma$ we obtain
approximations of the factors $m_{\gamma}$ which allow us to
compute the minimal polynomial $m_u$ and hence a geometric
solution of $\widehat{V}$. Our next result outlines this procedure
and estimates its complexity and error probability.

\begin{proposition}\label{prop: comput geo sol V hat}
Suppose that we are given a geometric solution of the variety
$V_{0,\gamma}$ for all $\gamma\in\Gamma$, as provided by Theorem
\ref{th: computation geo sol V_0}, with a linear form
$u\in\Q[\xon]$ whose coefficients are randomly chosen in the set
$\{1\klk 4\rho D^4\}$, where $\rho$ is a fixed positive integer.
Then there is an algorithm which computes a geometric solution of
the curve $\widehat{V}$ with error probability bounded by $1/\rho$
performing $O\big((n^2L+n^{1+\Omega}){\sf M}(\mathcal{M}_\Gamma)
{\sf M}(D)\big({\sf M}(D)+{\sf M}(E)\big)\big)$ arithmetic
operations in $\Q$. Here $L:=\max_{\gamma\in\Gamma}L_\gamma$,
where $L_\gamma$ is the number of arithmetic operations required
to evaluate the polynomials $h_{i,\gamma}$ of {\rm (\ref{eq: def h
gamma})} for all $\gamma\in\Gamma$ and
$\mathcal{M}_\Gamma:=\max_{\gamma\in\Gamma}\max\{\!\gamma_1\klk\!\gamma_{n+1}\}$.
\end{proposition}

\begin{proof} For each $\gamma\in \Gamma$, we apply the algorithm
underlying the proof of Proposition \ref{prop: approx minimal u
gamma} in order to obtain an approximation of $m_{u_\gamma}$ with
precision $2\gamma_{n+1}E$. Due to   (\ref{eq: minimal u gamma as
projection}), this polynomial immediately yields an approximation
with precision $2E$ of $m_{\gamma}(T,Y)$ in
${\Q}(\!(T^{1/\gamma_{n+1}})\!)[Y]$.

Multiplying all these approximations, we obtain an approximation
with precision $2E$ of the polynomial
$\widehat{m}_u=\prod_{\gamma\in\Gamma}m_{\gamma}$ of (\ref{eq:
minimal u, puiseux and cells}). Since every coefficient $a_j(T)$
of $\widehat{m}_u\in\Q(T)[Y]$ is a rational function of $\Q(T)$
having a reduced representation with numerator and denominator of
degree at most $E$, such a representation of $a_j(T)$ can be
computed from its approximation with precision $2E$ using Pad\'e
approximation with $O({\sf M}(E))$ arithmetic operations in $\Q$.

In order to estimate the complexity of the whole procedure, we
estimate the complexity of its three main steps:
\begin{enumerate}
    \item[$(i)$] the computation of the polynomials
    $m_{\gamma}$ with precision
    $2E$ for all $\gamma\in\Gamma$,
    which requires $O\big(\sum_{\gamma\in\Gamma}(nL_\gamma+n^{\Omega}){\sf
M}(D_\gamma)\big({\sf M}(M_\gamma){\sf
M}(D_\gamma)+E\gamma_{n+1}\big)\big)$ arithmetic operations in
$\Q$,
    \item[$(ii)$] the computation of the product $\prod_{\gamma\in\Gamma}m_{\gamma}$
    with precision $2E$, which requires $O\big({\sf M}(D){\sf M}(E)\big)$
    arithmetic operations in $\Q$,
    \item[$(iii)$] the computation of a reduced representation of all the
coefficients of $\widehat{m}_u\in\Q(T)[Y]$, which requires
$O\big({\sf M}(E)D\big)$ arithmetic operations in $\Q$.
\end{enumerate}
In conclusion, the algorithm performs $O\big((nL+n^\Omega) {\sf
M}(\mathcal{M}_\Gamma){\sf M}(D)\big({\sf M}(D)+{\sf
M}(E)\big)\big)$ arithmetic operations in $\Q$, where
$\mathcal{M}_\Gamma:=\max_{\gamma\in\Gamma}\{M_\gamma,\gamma_{n+1}\}$
and $L:=\max_{\gamma\in\Gamma}L_\gamma$.

Next we discuss how this procedure can be extended to the
computation of a geometric solution of $\widehat{V}$ in the sense
of Section \ref{subsect: geometric sol}. Two computations of the
above procedure involve divisions by the coefficients $u_i$ of the
linear form $u$: the computation of the resultant of (\ref{eq:
comput resultant m u tilde}) for all $\gamma\in\Gamma$ and the
Pad\'e approximations of $(iii)$. Both computations are reduced to
$D$ applications of the EEA, which is performed in a ring
$\Q(\Lambda)$. A similar analysis as in Proposition \ref{prop:
first step computation geo sol V_0} shows that all the
denominators in $\Q[\Lambda]$ arising during such application of
the EEA are divisors of a polynomial of degree $4D^4$. Therefore,
according to Lemma \ref{lemma: min/resgeom}, we conclude that a
geometric solution of $\widehat{V}$ can be computed with
$O\big((n^2L+n^{1+\Omega}){\sf M}(\mathcal{M}_\Gamma){\sf
M}(D)\big({\sf M}(D)+{\sf M}(E)\big)\big)$ arithmetic operations
in $\Q$, with an algorithm with error probability at most
$1/\rho$, provided that the coefficients of $u$ are randomly
chosen in the set $\{1\klk 4\rho D^4\}$.
\end{proof}

Putting together Theorem \ref{th: computation geo sol V_0} and
Proposition \ref{prop: comput geo sol V hat} we obtain the main
result of this section:
\begin{theorem}\label{th: comput geo sol V hat complete} Let
$\rho$ be a fixed positive integer. Suppose that the coefficients
of the linear form $\widetilde{u}$ of the statement of Theorem
\ref{th: computation geo sol V_0} and of the linear form $u$ are
randomly chosen in the set $\{1\klk 4n\rho D^4\}$. Then the
algorithm underlying Theorem \ref{th: computation geo sol V_0} and
Proposition \ref{prop: comput geo sol V hat} computes a geometric
solution of the curve $\widehat{V}$ with error probability
$3/\rho$ performing $O\big((n^2L+n^{1+\Omega}){\sf
M}(\mathcal{M}_\Gamma)\log(\mathcal{Q}){\sf M}(D)\big({\sf
M}(D)+{\sf M}(E)\big)\big)$ arithmetic operations in $\Q$. Here
$L:=\max_{\gamma\in\Gamma}L_\gamma$, where $L_\gamma$ is the
number of arithmetic operations required to evaluate the
polynomials $h_{i,\gamma}$ of {\rm (\ref{eq: def h gamma})} for
all $\gamma\in\Gamma$, $\mathcal{Q}:=2\max_{1\le i\le
n}\{\|q\|;q\in\Delta_i\}$, and
$\mathcal{M}_\Gamma:=\max_{\gamma\in\Gamma}\|\gamma\|$.
\end{theorem}
%
%
\subsection{Solving a sufficiently generic sparse system}
\label{subsect: solution generic sparse system}
Now we obtain a geometric solution of the zero-dimensional variety
$V_1:= \{ x \in \C^n: h_1(x) = 0, \dots, h_n(x) = 0\}$ from a
geometric solution of the curve $\widehat V$.

With notations as in the previous section, we have that $V_1 =
\pi^{-1}(1)$, where $\pi:\widehat V \to \A^1$ is the linear
projection defined by $\pi(x,t):=t$. Moreover, due to Lemma
\ref{lemma: V hat es curva}, the equality $V_1 = \pi^{-1}(1) \cap
\widehat V$ holds.

This enables us to easily obtain a geometric solution of $V_1$
from a geometric solution of the curve $\widehat V$. Indeed, let
$\widehat{m}_u(T,Y), \widehat{v}_1(T,Y)\klk \widehat{v}_n(T,Y)$ be
the polynomials which form a geometric solution of $\widehat V$
associated to a linear form $u\in\Q[X]$. Suppose further that the
linear form $u$ separates the points of $V_1$. Making the
substitution $T=1$, we obtain new polynomials
$\widehat{m}_u(1,Y),\widehat{v}_1(1,Y)\klk \widehat{v}_n(1,Y)\in
\Q[Y]$ such that $\widehat{m}_u(1,u(X))$ and $\frac{\partial
m_u}{\partial Y} (1, u(X)) X_i - \widehat{v}_i(1,u(X))$ $(1\le i
\le n)$ vanish over $V_1$. Taking into account that $\deg_Y(m_u) =
D  = \# V_1$ and that $u$ separates the points of $V_1$, it
follows that the polynomials
$\widehat{m}_u(1,Y),\widehat{v}_1(1,Y)\klk \widehat{v}_n(1,Y)\in
\Q[Y]$ form a geometric solution of $V_1$.
\begin{proposition}\label{prop: comput geo sol V 1}
Let $\rho$ be a fixed positive integer. With assumptions and
notations as in Theorem \ref{th: comput geo sol V hat complete},
the algorithm described above computes a geometric solution of the
zero-dimensional variety $V_1$ with error probability $4/\rho$
using $O\big((n^2L+n^{1+\Omega}){\sf
M}(\mathcal{M}_\Gamma)\log(\mathcal{Q}){\sf M}(D)\big({\sf
M}(D)+{\sf M}(E)\big)\big)$ arithmetic operations in $\Q$.
\end{proposition}
%
%
\section{The solution of the original system}
\label{sect: solution original system}
Let notations and assumptions be as in the previous sections.
Assume that we are given a geometric solution $m_u(Y), v_1(Y),
\dots, v_n(Y)$ of the zero-dimensional variety $V_1$ defined by
the polynomials $h_1:=f_1 +g_1,\dots, h_n:= f_n + g_n$. Assume
further that the linear form $u$ of such a geometric solution
separates the points of the zero--dimensional variety
$f_1=\cdots=f_n=0$. In this section we describe a procedure for
computing a geometric solution of the input system
$f_1=\dots=f_n=0$.

For this purpose, we introduce an indeterminate $T$ over $\Q[X]$
and consider the ``deformation'' $F_1\klk F_n\in\Q[X,T]$ of the
polynomials $f_1\klk f_n$ defined in the following way:
$$F_i(X, T):=f_i(X)+(1-T)g_i(X)\quad (1\le i\le n).$$
Set $\mathcal{V}:=\{(x,t)\in\A^{n+1}:F_1(x,t)=\dots=F_n(x,t)=0\}$
and denote by $\pi:\mathcal{V}\to\A^1$ the projection map defined
by $\pi(x,t):=t$. As in Subsection \ref{subsect: comput geo sol V
hat}, we introduce the variety $\mathcal{V}_{\rm dom} \subset
\A^{n+1}$ defined as the union of all the irreducible components
of $\mathcal{V}$ whose projection over $\A^1$ is dominant.
%
%
\subsection{Solution of the second deformation.} In this section
we describe an efficient procedure for computing a geometric
solution of $\mathcal{V}_{\rm dom}$ from the geometric solution of
$\pi^{-1}(0)$ provided by Proposition \ref{prop: comput geo sol V
1}.

Since $\pi^{-1}(0)$ is the variety defined by the ``sufficiently
generic'' sparse system $h_1(X)=F_1(X, 0)=0,\dots,
h_n(X)=F_n(X,0)=0$, with similar arguments to those leading to the
proof of Lemma \ref{lemma: V hat es curva}, it is not difficult to
see that the polynomials $F_1,\dots, F_n$, the variety
$\mathcal{V}$, the projection $\pi: \mathcal{V} \to \A^1$, and the
fiber $\pi^{-1}(0)$ satisfy all the assumptions of Lemma
\ref{lemma: variedad dom}. We conclude that $\mathcal{V}_{\rm
dom}$ is a curve and  that the identity $\mathcal{V} \cap
\pi^{-1}(0) = \mathcal{V}_{\rm dom}\cap \pi^{-1}(0)$ holds.
Furthermore, Lemma \ref{lemma: variedad dom} implies that all the
hypotheses of \cite[Theorem 2]{Schost03} are satisfied.

Therefore, applying the ``formal Newton lifting process''
underlying the proof of \cite[Theorem 2]{Schost03}, we compute
polynomials $m(T,Y),v_1(T,Y)\klk v_n(T,Y)\in\Q[T,Y]$ which form a
geometric solution of $\mathcal{V}_{\rm dom}$. The formal Newton
lifting process requires $O\big((n{L}^\prime+n^{\Omega+1}){\sf
M}(D){\sf M}(E^\prime)\big)$ arithmetic operations in ${\Q}$,
where ${L}^\prime$ denotes the number of arithmetic operations
required to evaluate $F_1\klk F_n$ and $E^\prime$ is any upper
bound of the degree of $m$ in the variable $T$.

We can apply Lemma \ref{lemma: estimate sparse height} in order to
estimate $\deg_T m$ in combinatorial terms. Indeed, let
$\widetilde{Q}_1,\dots, \widetilde{Q}_n\subset \R^{n+1}$ be the
Newton polytopes of the polynomials $F_1,\dots,F_n$ and let
$\Delta\subset \R^{n+1}$ be the standard unitary simplex in the
plane $\{T = 0\}$. Since $\widetilde{Q}_i\subset Q_i\times[0,1]$
holds for $1\le i\le n$, where $Q_i\subset\R^n$ is the Newton
polytope of $h_i$, by (\ref{eq: estimate height second def}) of
Lemma \ref{lemma: estimate sparse height} we deduce the following
estimate:
\begin{equation}\label{eq: estimate height second def bis}
\deg_T m_u\le E':=\sum_{i=1}^nMV(\Delta,Q_1,\dots,Q_{i-1},Q_{i+1}
\klk Q_n).
\end{equation}
With this definition of $E'$, we have:
\begin{proposition}\label{prop: comput geo sol V dom}
Suppose that we are given a geometric solution of the variety
$V_1$, as provided by Proposition \ref{prop: comput geo sol V 1}.
A geometric solution of $\mathcal{V}_{\rm dom}$ can be
deterministically computed with
$O\big((nL^\prime+n^{\Omega+1}){\sf M}(D){\sf M}(E^\prime)\big)$
arithmetic operations in $\Q$.
\end{proposition}
%
%
\subsection{Solving the input system.} Making the
substitution $T=1$ in the polynomials $m(T,Y),v_i(T,Y)$ $(1\le
i\le n)$ which form the geometric solution of $\mathcal{V}_{\rm
dom}$ computed by the algorithm of Proposition \ref{prop: comput
geo sol V dom} we obtain polynomials $m(1,Y),v_1(1,Y),$ $\dots,
v_n(1,Y)\in\Q[Y]$ which represent a complete description of our
input system $f_1(X)=\cdots=f_n(X)=0$, eventually including
multiplicities. Such multiplicities are represented by multiple
factors of $m(1,Y)$, which are also factors of $v_1(1,Y)\klk
v_n(1,Y)$ (see e.g. \cite[\S 6.5]{GiLeSa01}). In order to remove
them, we compute $a(Y):=\mathrm{gcd}\big(m(1,Y),(\partial
m/\partial Y)(1,Y)\big)$, and the polynomials $m(1,Y)/a(Y)$,
$(\partial m/\partial Y)(1,Y)/a(Y)$, $v_i(1,Y)/a(Y)$ $(1\le i\le
n)$. These polynomials form a geometric solution of our input
system and can be computed with $O\big(n{\sf M}(D)E'\big)$
additional arithmetic operations in $\Q$.

Summarizing, we sketch the whole procedure computing a geometric
solution of the input system $f_1=\dots=f_n=0$. Fix $\rho\ge 4$.
We randomly choose the coefficients of the polynomials $g_1\klk
g_n$ in the set $\{1\klk 4\rho(nd)^{2n+1}+2\rho
n^22^{\mathcal{N}_1\plp \mathcal{N}_s}\}$ and coefficients of
linear forms $u,\widetilde{u}$ in the set $\{1\klk 16n\rho D^4\}$.
By Theorem \ref{th: Zippel-Schwartz} it follows that the
polynomials $g_1\klk g_n$ and the linear forms $u,\widetilde{u}$
satisfy all the conditions required with probability at least
$1-1/\rho$. Then we apply the algorithms underlying Propositions
\ref{prop: comput geo sol V 1} and \ref{prop: comput geo sol V
dom} in order to obtain a geometric solution of the variety
$\mathcal{V}_{\rm dom}$. Finally, we use the procedure above to
compute a geometric solution of the input system
$f_1=\dots=f_n=0$. This yields the following result:
\begin{theorem}\label{th: comput_geo_sol input system}
The algorithm sketched above computes a geometric solution of the
input system $f_1=\dots=f_n=0$ with error probability at most
$1/\rho$ using
$$O\Big(\big(n^2\max\{L,L'\}+n^{1+\Omega}\big)\,{\sf M}(D)
\,\big(\log(\mathcal{Q}){\sf M}(\mathcal{M}_\Gamma)\big({\sf
M}(D)+{\sf M}(E)\big)+{\sf M}(E')\big)\Big)$$
arithmetic operations in $\Q$. Here
$L:=\max_{\gamma\in\Gamma}L_\gamma$, where $L_\gamma$ is the
number of arithmetic operations required to evaluate the
polynomials $h_{i,\gamma}$ of {\rm (\ref{eq: def h gamma})} for
all $\gamma\in\Gamma$, ${L}^\prime$ denotes the number of
arithmetic operations required to evaluate $F_1\klk F_n$,
$\mathcal{M}_\Gamma:=\max_{\gamma\in\Gamma}\|\gamma\|$, and
$\mathcal{Q}:=2\max_{1\le i\le n}\{\|q\|;q\in\Delta_i\}$.
\end{theorem}
%
%
\providecommand{\bysame}{\leavevmode\hbox
to3em{\hrulefill}\thinspace}
\providecommand{\MR}{\relax\ifhmode\unskip\space\fi MR }
\providecommand{\MRhref}[2]{%
  \href{http://www.ams.org/mathscinet-getitem?mr=#1}{#2}
} \providecommand{\href}[2]{#2}

\end{document}